\documentclass[11pt]{article}

\usepackage[utf8]{inputenc}

\usepackage{times,amssymb,amsmath,exscale,array,latexsym}
\usepackage{graphicx}
\usepackage{epsfig}
\usepackage{color}
\usepackage[dvipsnames]{xcolor}
\usepackage{subfig}
\definecolor{marin}{rgb}   {0.,   0.3,   0.7}
\definecolor{rouge}{rgb}   {0.8,   0.,   0.}
\definecolor{sepia}{rgb}   {0.8,   0.5,   0.}
\usepackage[colorlinks,citecolor=marin,linkcolor=rouge,
            bookmarksopen,
            bookmarksnumbered
           ]{hyperref}

\newcommand{\e}{\ensuremath{\mathrm{e}}}
\newcommand{\tr}{\ensuremath{\mathrm{tr}}}

\newcommand{\psis}{{\mathcal S}}
\renewcommand{\r}{\mathbb{R}}

\newcommand{\psiba}{\chi^*}
\newcommand{\psib}{\chi}

\numberwithin{equation}{section}

\newcommand{\QED}{\mbox{}\hfill \raisebox{-0.2pt}{\rule{5.6pt}{6pt}\rule{0pt}{0pt}}
          \medskip\par}

\newcommand{\cL}{{\mathcal L}}

\newcommand{\blau}[1]{{\color{MidnightBlue} #1}}
\makeindex
\begin{document}
\title{Families of efficient low order processed composition methods
}

\author{S. Blanes$^{1}$, F. Casas$^{2}$, A. Escorihuela-Tom\`as$^{3}$ \\[2ex]
$^{1}$ {\small\it Universitat Polit\`ecnica de Val\`encia, Instituto de Matem\'atica Multidisciplinar, 46022-Valencia, Spain}\\{
\small\it email: serblaza@imm.upv.es}\\[1ex]
$^{2}$ {\small\it Departament de Matem\`atiques and IMAC, Universitat Jaume I, 12071-Castell\'on, Spain}\\{
\small\it email: Fernando.Casas@mat.uji.es}\\[1ex]
$^{3}$ {\small\it Departament de Matem\`atiques, Universitat Jaume I, 12071-Castell\'on, Spain}\\{
\small\it email: alescori@uji.es}\\[1ex]
}


%
\maketitle

\begin{abstract}

New families of composition methods with processing of order 4 and 6 are presented and analyzed. They are specifically designed to be used for
the numerical integration of differential equations whose vector field is separated into three or more parts which are explicitly solvable. The new schemes
are shown to be more efficient than previous state-of-the-art splitting methods.



\end{abstract}\bigskip

\noindent \emph{AMS numbers}: 65L05, 65P10

\noindent \emph{Keywords}: composition methods, effective error, processing technique, splitting, Hamiltonian systems

\section{Introduction}
\label{sec.1}

Structure-preserving numerical integration methods are nowadays a common tool in many areas of physics, chemistry and computational mathematics
\cite{blanes16aci,hairer06gni}. 
Among them, splitting methods constitute a natural option when the differential system can be separated into two parts, so that
each of them is explicitly integrable \cite{blanes24smf,mclachlan02sm}. Suppose the vector field $f$ in 
\begin{equation}   \label{ode.1}
   \dot{x} \equiv \frac{dx}{dt} = f(x),  \qquad x(t_0) = x_0 \in \mathbb{R}^D
 \end{equation}
can be split into two parts, $f(x) = f_1(x) + f_2(x)$, so that each subproblem
\[
   \dot{x} = f_i(x), \qquad x(t_0) = x_0, \qquad i=1, 2
\]
is explicitly solvable, with solution $x(t) = \varphi_t^{[i]}(x_0)$. Then, the composition
\begin{equation} \label{chi2}
  \chi_h = \varphi_h^{[2]} \circ \varphi_h^{[1]} 
\end{equation}
provides a first-order approximation to the exact flow $\varphi_h$ of (\ref{ode.1}):
\begin{equation} \label{lt}
  \chi_h(x_0) =   \varphi_h(x_0) + \mathcal{O}(h^2), \qquad \mbox{ as } \quad h \rightarrow 0,
\end{equation}
whereas the palindromic composition
\begin{equation} \label{lf}
  \mathcal{S}_h^{[2]} = \varphi_{h/2}^{[1]} \circ  \varphi_{h}^{[2]} \circ  \varphi_{h/2}^{[1]}, 
\end{equation}
known as the Strang splitting method, is of second-order. Higher-order schemes can be constructed as 
\[
  \psi_h =  \varphi_{h b_s}^{[2]} \circ \varphi_{h a_s}^{[1]}  \circ \cdots \circ   \varphi_{h b_1}^{[2]} \circ \varphi_{h a_1}^{[1]}
\]
if the coefficients $a_j$, $b_j$ are conveniently chosen so as to satisfy the required order conditions \cite{hairer06gni,mclachlan02sm}.

There are problems, however,
where $f$ has to be split into more than two terms for each part to be explicitly integrable, i.e., 
$f(x) =  \sum_{i=1}^n f_i(x)$, $n \ge 3$. In that case, method \eqref{lt} generalizes to
\begin{equation}  \label{LT}
   \chi_h = \varphi^{[n]}_{h}\circ \, \varphi^{[n-1]}_{h}\circ \, \cdots \, \circ \varphi^{[2]}_{h} \circ \varphi^{[1]}_{h} 
\end{equation}
(or any other permutation of the sub-flows $\varphi_h^{[j]}$), leading to a first-order approximation.
This is also the case of the adjoint of $\chi_h$, defined as $\chi_h^* = (\chi_{-h})^{-1}$, namely
\[
   \chi_h^* = \varphi^{[1]}_{h}\circ \, \varphi^{[2]}_{h}\circ \, \cdots \, \circ \varphi^{[n-1]}_{h} \circ \varphi^{[n]}_{h}, 
\]
whereas the composition $\chi_{h/2} \circ \chi_{h/2}^*$ leads to a second-order approximation $\mathcal{S}_h^{[2]}$, 
the generalization of scheme (\ref{lf}) to this setting. Higher order integrators can be constructed as compositions of
$\mathcal{S}_h^{[2]}$. Thus, in particular, the
4th-order scheme
\begin{equation} \label{yoshida}
  \mathcal{S}_h^{[4]} = \mathcal{S}_{\alpha_1 h}^{[2]} \circ \mathcal{S}_{\alpha_2 h}^{[2]} \circ \mathcal{S}_{\alpha_1 h}^{[2]}, \qquad
  \mbox{ with } \qquad  \alpha_1 = \frac{1}{2 - 2^{1/3}}, \qquad \alpha_2 = 1 - 2 \alpha_1,
\end{equation}
is very popular in applications \cite{yoshida90coh,blanes08sac}, although it has large truncation errors and a short stability interval. 
Alternatively, the 4th-order scheme
\[
 \psi_h = \psis_{\alpha_1 h}^{[2]} \circ  \psis_{\alpha_1 h}^{[2]} \circ   \psis_{\alpha_3 h}^{[2]} \circ    \psis_{\alpha_1 h}^{[2]} \circ \psis_{\alpha_1 h}^{[2]},
  \qquad \mbox{ with } \qquad \alpha_1 = \frac{1}{4-4^{1/3}},  \quad \alpha_3 = 1 - 4 \alpha_1,
\]
first proposed in \cite{suzuki90fdo} and analyzed in detail in \cite{mclachlan02foh}, is much more efficient than (\ref{yoshida}), even if it
requires more computational effort per step, whereas 
palindromic compositions of the form
\begin{equation} \label{eq.2.1.2}
  \psi_h = \psis_{\alpha_{m} h}^{[2]}\circ \,  \psis_{\alpha_{m-1} h}^{[2]} \circ
\cdots\circ
 \psis_{\alpha_{2} h}^{[2]}  \circ \, \psis_{\alpha_{1}h}^{[2]}  \qquad \mbox{ with } \qquad 
(\alpha_1,\ldots,\alpha_{m}) \in \r^{m}
\end{equation}
and $\alpha_{m+1-i} = \alpha_i$, lead to efficient integrators of order $r > 6$ \cite{hairer06gni}. 

The most general situation corresponds to integrators of the form
\begin{equation} \label{eq.2.1.1}
\psi_h = \psib_{\alpha_{2s} h}\circ
\psiba_{\alpha_{2s-1}h}\circ
\cdots\circ
\psib_{\alpha_{2}h}\circ
\psiba_{\alpha_{1}h}, \quad \mbox{ with } \quad (\alpha_1,\ldots,\alpha_{2s})
   \in \r^{2s}
\end{equation}
verifying in addition the condition $\alpha_{2s+1-i} = \alpha_i$ to preserve time-symmetry. In fact, methods (\ref{eq.2.1.2}) constitute a particular case
of (\ref{eq.2.1.1}) when $\mathcal{S}_h^{[2]}$ is the Strang map (\ref{lf}). Among the most efficient 4th- and 6th-order methods of this class we can mention the schemes  introduced by Blanes and Moan 
\cite{blanes02psp}, denoted here as BM${_6^{[4]}}$ and BM$_{10}^{[6]}$, respectively. 
They have $s=6$ and $s=10$ \emph{stages}, and their coefficients are collected in Table \ref{table.1}.
The main truncation error of the 4th-order scheme BM${_6^{[4]}}$ is around 500 times smaller than
the error of method (\ref{yoshida}), thus compensating its higher computational cost per step. For future reference, this scheme reads  explicitly
\begin{equation}  \label{m.6}
 \mbox{ BM${_6^{[4]}}$ }: \quad  \chi_{\alpha_1 h} \circ \chi_{\alpha_2 h}^* \circ \chi_{\alpha_3 h} \circ \chi_{\alpha_4 h}^* \circ \chi_{\alpha_5 h} \circ 
  \chi_{\alpha_6 h}^* \circ \chi_{\alpha_6 h} \circ
  \chi_{\alpha_5 h}^* \circ
  \chi_{\alpha_4 h} \circ  
  \chi_{\alpha_3 h}^* \circ \chi_{\alpha_2 h} \circ \chi_{\alpha_1 h}^*.
\end{equation}

\begin{table}[t]
  \centering
    \renewcommand\arraystretch{1.1}
    \begin{tabular}{lll}
      \multicolumn{3}{c}{BM$_{6}^{[4]}$ (order 4)}\\
      \hline 
      $\alpha_1= 0.0792036964311957$ &\qquad $\alpha_2= 0.1303114101821663$ \\
      $\alpha_3= 0.22286149586760773$ & \qquad  $\alpha_4= -0.36671326904742574$ \\
      $\alpha_5=0.32464818868970624$ & \qquad $\alpha_6=0.10968847787674973$ \\
      \hline 
                                    & & \\
   \multicolumn{3}{c}{BM$_{10}^{[6]}$ (order 6)}\\
    \hline
    $\displaystyle \alpha_1= 0.0502627644003922$ &\qquad $\displaystyle \alpha_2= 0.0985536835006498$ \\
    $\displaystyle \alpha_3= 0.31496061692769417$ & \qquad $\displaystyle \alpha_4= -0.44734648269547816$ \\
    $\displaystyle \alpha_5=0.49242637248987586$ & \qquad $\displaystyle \alpha_6=-0.42511876779769087$ \\
    $\alpha_7 = 0.23706391397812188$  & \qquad $\alpha_8 = 0.19560248860005314$  \\
    $\alpha_9 = 0.34635818985072686$ & \qquad $\alpha_{10} = -0.36276277925434486$ \\
    \hline
    \end{tabular}
  \caption{\small{Fourth- and sixth-order symmetric (palindromic)
  composition methods of the form (\ref{eq.2.1.1}). They correspond to the splitting schemes \textit{S}$_6$ and
  \textit{S}$_{10}$ of \cite[Table 2]{blanes02psp}, respectively. \label{table.1}}}
\end{table}

Actually, schemes BM${_6^{[4]}}$ and BM$_{10}^{[6]}$ are constructed in \cite{blanes02psp} as splitting methods when $f$ is separated into two parts. 
Specifically, if
$\chi_h$ is taken as \eqref{chi2}, then scheme \eqref{eq.2.1.1} can be rewritten as 
\begin{equation} \label{eq:splitting}
\psi_h = \varphi^{[1]}_{a_{s+1} h}\circ \varphi^{[2]}_{b_{s}h}\circ \varphi^{[1]}_{a_{s} h}\circ \cdots\circ
 \varphi^{[1]}_{a_{2}h}\circ
 \varphi^{[2]}_{b_{1}h} \circ \varphi^{[1]}_{a_{1}h},
\end{equation}
where $a_{1}=\alpha_{1}$, and for $j=1,\ldots,s$,
\begin{equation}
  \label{eq:abalpha}
a_{j+1} = \alpha_{2j} + \alpha_{2j+1}, \qquad\quad  b_j=\alpha_{2j-1}+\alpha_{2j}
\end{equation}
 (with $\alpha_{2s+1}=0$). Conversely, any integrator of the form (\ref{eq:splitting}) satisfying the condition
$\sum_{j=1}^{s+1} a_j = \sum_{j=1}^{s} b_j$
  can be expressed in the form (\ref{eq.2.1.1}), as shown
  in \cite{mclachlan95otn}.

The aim of this paper is to 
construct new classes of 4th- and 6th-order composition methods that are even more efficient for problems that can be expressed
as the sum of three or more explicitly integrable terms. These are built by applying the processing technique.
Although a detailed study of composition methods with processing was carried out in \cite{blanes06cmf} and methods up
to order 12 were presented there, we think methods of low order require an additional treatment to improve their overall efficiency
due to their relevance in practical applications. This is the subject of the present work, where we present new families of processed
composition schemes with a better performance than the \emph{de facto} state-of-the-art numerical integrators of Table \ref{table.1}. 
This claim is further substantiated by  numerical tests on different problems. Although the new schemes are specifically designed
and optimized for systems separable into $n \ge 3$ parts, the case $n=2$ is also contemplated.
  Furthermore, the new schemes are also compared with the most efficient processed methods of order 4 and 6 obtained in \cite{blanes06cmf}, whose
  coefficients (for the kernel) are collected in Table \ref{table.1.2}.
  
\begin{table}[t]
  \centering
    \renewcommand\arraystretch{1.1}
    \begin{tabular}{lll}
      \multicolumn{3}{c}{BCM$_{6}^{[4]}$ (effective order 4)}\\
      \hline 
      $\alpha_1= \alpha_2 =\alpha_3=\alpha_4$ &\qquad $\alpha_4= 0.1341940158122142$ \\
      $\alpha_5= -0.3141940158122142$ & \qquad  $\alpha_6= 0.27741795256335733$ \\
      \hline 
                                    & & \\
   \multicolumn{3}{c}{BCM$_{9}^{[6]}$ (effective order 6)}\\
    \hline
    $\displaystyle \alpha_1= \alpha_2=\alpha_3=\alpha_4=\alpha_5$ &\qquad $\displaystyle \alpha_5= 0.1106570871853300$ \\
    $\displaystyle \alpha_6= -0.2854111127287940$ & \qquad $\displaystyle \alpha_7= 0.2138498496192465$ \\
    $\displaystyle \alpha_8=-0.3402583791791715$ & \qquad $\displaystyle \alpha_9=0.35853420636206895$ \\
    \hline
    \end{tabular}
    \caption{\small{Schemes of effective order four and six of the form (\ref{eq.2.1.1}).  They correspond to methods $P_64$ and $P_96$
        of \cite[Table 5]{blanes06cmf}, respectively. With an appropriate processor, they render 4th- and 6th-order composition schemes. \label{table.1.2}}}
\end{table}

\section{Processed composition methods}

Processed (or corrected) methods are of the form
\begin{equation} \label{proces.1}
  \hat{\psi}_h = \pi_h \circ \psi_h \circ \pi_h^{-1}.
\end{equation}
The integrator $\psi_h$ is called the \emph{kernel} and the (near-identity) map $\pi_h$ is the \emph{processor} or \emph{corrector}.
The method $\psi_h$ is said to be of effective order $r$ if a processor exists such that $\hat{\psi}_h$ is of (conventional) order $r$ \cite{butcher96tno}. Note that, since
\[
  \hat{\psi}_h^N = \pi_h \circ \psi_h^N \circ \pi_h^{-1},
\]
applying  $\hat{\psi}_h$ over $N$ steps with constant $h$ only involves $N$ evaluations of the kernel, whereas $\pi_h^{-1}$ is
computed only at the beginning and $\pi_h$ when output is desired \cite{lopezmarcos97esi,blanes99siw}.

Processed integrators have shown to be very efficient in a variety of systems, ranging from
near-integrable problems to Hamiltonian systems separable into kinetic and potential energy \cite{lopezmarcos97esi,wisdom96sco,blanes99siw}. This is due essentially to the fact that the kernel
has to satisfy a much reduced set of order conditions (in other words, some of the order conditions can be fulfilled by the processor) and therefore
they require less computational effort than a conventional method of the same order.

The derivation and analysis of the effective order conditions for kernels of the form (\ref{eq.2.1.1}) has been done in \cite{blanes06cmf},
where kernels of effective orders 4 and 6 have also been proposed (for order $r > 6$ it is more advantageous to consider directly palindromic
compositions of the form (\ref{eq.2.1.2})). Here we briefly summarize the treatment when $r \le 6$ and construct new families of
schemes.


The corresponding analysis can be carried out with the help of the Lie formalism. To proceed, we introduce the Lie derivative $F$ associated with
$f$ in (\ref{ode.1}), and defined as
\[
    F \, g (x) = \left.\frac{d}{d h}\right|_{h=0} g(\varphi_{h}(x)) 
\]
for each smooth function $g:\mathbb{R}^D \rightarrow \mathbb{R}$ and $x \in \mathbb{R}^D$, that is,
\begin{equation}   \label{eq:3b}
F\, g (x) = f(x) \cdot \nabla g(x).
\end{equation}
Then, the $h$-flow of Eq.  (\ref{ode.1}) verifies \cite{hairer06gni,sanz-serna94nhp}
\[
   g(\varphi_{h}(x)) = \e^{h F}g (x),
\]   
where $\e^{h F}$ is defined as a series of linear differential operators
\[
   \e^{h F} = \sum_{k=0}^{\infty} \frac{h^k}{k!} F^k.
\]
Analogously, for the basic method $\chi_h$ of (\ref{LT}), one can associate a series of linear operators so that \cite{blanes08sac}
\[
    g(\psib_h(x)) = \e^{Y(h)}g(x), \qquad \mbox{ with }  \qquad Y(h) = \sum_{k \ge 1} h^k Y_k
\]
with $Y_1 = F$, whereas for its adjoint one has $g(\psiba_h(x)) = \e^{-Y(-h)} g(x)$.
In consequence, the integrator (\ref{eq.2.1.1}) has associated a series 
 $\Psi(h)$   of differential operators given by
\begin{equation} \label{series}
\Psi(h) = \e^{-Y(-h\alpha_1)} \, \e^{Y(h \alpha_2)} \cdots \, 
\e^{-Y(-h\alpha_{2s-1})} \,\e^{Y( h \alpha_{2s})},
\end{equation}
in the sense that $g(\psi_{h}(x)) = \Psi(h)\,  g (x)$.
Successive applications of the Baker--Campbell--Hausdorff formula \cite{varadarajan84lgl} in (\ref{series}) allow us to formally express $\Psi(h)$ as
only one exponential,
\begin{equation}   \label{modH}
   \Psi(h) = \exp(K(h)), \qquad \mbox{ with } \qquad K(h) = \sum_{j \ge 1} h^j K_j,
\end{equation}   
the terms $h^j K_j \in \mathcal{L}_j$ for each $j\geq 1$
and
$\mathcal{L} = \bigoplus_{j\geq 1} \mathcal{L}_j$ is the graded free Lie algebra generated by $\{Y_1, Y_3, Y_5, \ldots \}$ \cite{munthe-kaas99cia}.
For the particular basis of $ \mathcal{L}_j$, $j \le 7$,  collected in Table \ref{table.2}, 
the operator $K(h)$ associated with a consistent time-symmetric method (\ref{eq.2.1.1}) reads
\begin{equation} \label{expreK}
   K(h)  =  h F + h^3 \big( k_{3,1} E_{3,1} + k_{3,2} E_{3,2} \big) + h^5 \sum_{\ell=1}^6 k_{5,\ell}  E_{5,\ell} 
      + h^7 \sum_{\ell=1}^{18} k_{7,\ell}  E_{7,\ell} + \mathcal{O}(h^9),
\end{equation}  
where $k_{n, \ell}$ are polynomials of degree $n$ in the coefficients $\alpha_j$. Due to the time-symmetry, only odd order terms appear in \eqref{expreK}
\cite{sanz-serna94nhp}.


\begin{table}[h]
\small
\begin{center}
  \begin{tabular}{c|lll}
    $\cL_j$ & \multicolumn{3}{c}{{\rm Basis of} $\cL_j$}   \\    \hline
$\cL_1$ & $E_{1,1}=Y_1$ &   &      \\   \hline
$\cL_2$ & $E_{2,1}=Y_2$ &   &      \\   \hline
$\cL_3$ & $E_{3,1}=Y_3$   & $E_{3,2}=[Y_1, Y_2]$  &      \\   \hline
$\cL_4$ & $E_{4,1} = Y_4$  & $E_{4,2}=[Y_1, Y_3]$  &  $E_{4,3}=[Y_1, E_{3,2}]$      \\   \hline
$\cL_5$ & $E_{5,1}=Y_5$  & $E_{5,2}=[Y_1,Y_4]$ & $E_{5,3}=[Y_1,E_{4,2}]$   \\ 
              & $E_{5,4}=[Y_1, E_{4,3}]$ & $E_{5,5}=[Y_2, Y_3]$ & $E_{5,6}=[Y_2, E_{3,2}]$ \\ \hline
    $\cL_6$ & $E_{6,1}=Y_6$  & $E_{6,2}=[Y_1,Y_5]$  &   \\  
             & $E_{6,3+i}=[Y_1,E_{5,2+i}]$&  $i=0,\ldots,4$ & \\
             & $E_{6,8}=[Y_2,Y_4]$  & $E_{6,9}=[Y_2,E_{4,2}]$  &     \\  \hline
    $\cL_7$ & $E_{7,1}=Y_7$ & $E_{7,2}=[Y_1,Y_6]$ & \\
             & $E_{7,3+j}=[Y_1,E_{6,2+j}]$  & $j=0,\ldots,7$  &     \\ 
             & $E_{7,11}=[Y_2,Y_5]$& $E_{7,12+k} = [Y_2, E_{5,2+k}]$ & $k=0,\ldots,4$   \\
             & $E_{7,17}=[Y_3,Y_4]$ &  $E_{7,18}=[Y_3,E_{4,2}]$ & \\ 
\end{tabular}
\caption{Particular basis of $\cL_j$, $1 \le j \le 7$, taken in this work. }\label{table.2}
\end{center}
\end{table}

 As for the processor, when the kernel is time-symmetric, it can be chosen such that $\pi_{-h} = \pi_h$. In that case, its associated series of
 differential operators can be expressed as $\Pi(h) = \exp(P(h))$, with $P(-h) = P(h)$ \cite{blanes04otn}. Explicitly, 
\begin{equation} \label{p.series}
   P(h)  =   h^2 p_{2,1}  E_{2,1} + h^4 \big( p_{4,1} E_{4,1} + p_{4,2} E_{4,2} 
   + p_{4,3} E_{4,3} \big) +  \mathcal{O}(h^6).
\end{equation}
In consequence, the series of operators associated to the processed method (\ref{proces.1}) is
\begin{equation} \label{cm.series}
  \hat{\Psi}(h) = \e^{\hat{F}(h)} = \Pi(h)^{-1} \, \Psi(h) \, \Pi(h) = \e^{-P(h)} \, \e^{K(h)} \, \e^{P(h)},
\end{equation}

so that
\begin{equation} \label{f.hat}
\begin{aligned}
  & \hat{F}(h) = \e^{-\mathrm{ad}_{P(h)}} K(h) = K(h) - [P(h), K(h)] + \frac{1}{2} [P(h),[P(h), K(h)]] + \cdots \\
  & \quad = F +  h^3 \big( f_{3,1} E_{3,1} + f_{3,2} E_{3,2} \big) + h^5 \sum_{\ell=1}^6 f_{5,\ell}  E_{5,\ell} + \mathcal{O}(h^7),
\end{aligned}  
\end{equation}
where $f_{i,j}$ depend on $k_{\ell,m}$ and $p_{k,n}$. The method is of order $r$ if $f_{i,j} = 0$ for all $1 < i \le r$. This leads
to a set of restrictions on the coefficients of the kernel (the effective order conditions), whereas $P(h)$ is obtained by fixing its
coefficients $p_{k,n}$ so as to satisfy the remaining conditions. Since the kernel is time-symmetric, we can choose without
loss of generality $p_{2n-1,j} = 0$, $n \ge 1$  \cite{blanes99siw}.
All these conditions up to order 7 are collected in Table 
\ref{table.3}.

\begin{table}[htp]
\small
\begin{center}
\begin{tabular}{l}
  \begin{tabular}{c|lll}
 $r$ & {\rm Eff. order conditions} & \; {\rm Processor conditions} &  \\
\hline
 1  & $k_{1,1}=1$ &   &  \\
 3  & $k_{3,1}=0$ &    \;  $p_{2,1} = k_{3,2}$   & \\
 5  & $k_{5,1}=k_{5,5}=0$  &   &    \\
     & $k_{5,6} = - \frac{1}{2} k_{3,2}^2$ &  \; $p_{4,1+i} =  k_{5,2+i}, \qquad i= 0,1,2$  &\\
 7  & $k_{7,i}=0, \quad i = 1,11,15,17,18$ &\; $p_{6,1+l} = k_{7,2+l}, \qquad l= 0,\ldots,5$   &\\
     & $k_{7,16} = \frac{1}{6} k_{3,2}^3$ &  \; $p_{6,7} = k_{7,8}-\frac{1}{2}k_{3,2}k_{5,4}$ &\\
     & $k_{7,10+j} = -k_{3,2} k_{5,j}, \quad j = 2,3,4$ &  \; $p_{6,8} = k_{7,9}-\frac{1}{2}k_{3,2}k_{5,2}$&  \; $p_{6,9} = k_{7,10}-\frac{1}{2}k_{3,2}k_{5,3}$ \\
\end{tabular}
\end{tabular}
\caption{Order conditions for the kernel $\psi_h$ and processor $\pi_h$ up to order $r$ when the time-symmetric kernel ($k_{2n,i}=0$) is
expressed in the basis of Table \ref{table.2}. In that case, $p_{2n-1,j}=0$, $n=1,2,3$. } \label{table.3}
\end{center}
\end{table}

\section{Processed methods of order 4}
\label{sect3}

\subsection{Construction of kernels}

According with Table \ref{table.3}, two order conditions have to be solved by a palindromic kernel of the form (\ref{eq.2.1.1}) to achieve effective order 4. 
In other words, $s \ge 2$. As noticed previously \cite{blanes06cmf}, with $s=2$ there are only complex solutions, so that more stages
have to be introduced and, consequently, one has free parameters. In these circumstances, 
 some criterion has to be adopted to make an appropriate selection of the free parameters. 
 A standard strategy consists in minimizing the non-correctable error terms at order
5, since they cannot be removed by a processor anyway. In our case, these terms can be grouped into the function
\begin{equation} \label{error5}
   \mathcal{E}_5 \equiv \left(k_{5,1}^2 + k_{5,5}^2 + (k_{5,6}+ \frac{1}{2}k_{3,2}^2)^2 \right)^{1/2},
\end{equation}
and also corresponds to the 2-norm of the vector   $(f_{5,1},\ldots, f_{5,6})$ of the coefficients $f_{5,j}$
in the series (\ref{f.hat}),
once the coefficients $p_{i,j}$ are chosen according with the prescription of Table \ref{table.3} up to $i=4$. 
To take into account the computational cost corresponding to a kernel (\ref{eq.2.1.1}) with $s$ stages, the following
function, with $r=5$, is used to measure the relative efficiency of each method \cite{blanes06cmf}:
\begin{equation} \label{effic}
  E_{\mathrm{ef}}^{(r)} = s \, \mathcal{E}_r^{1/(r-1)}.
\end{equation}
In addition, we also keep track of the non-correctable error terms at orders 7
{and 9 with the same function (with $r=7$ and $r=9$),
and the size of the coefficients, as measured by the 1-norm of the vector $(\alpha_1, \ldots, \alpha_{2s})$. In this way we hope
to keep higher-order errors under control.

We have analyzed compositions with $3 \le s \le 9$ stages. Specifically, 
we have solved the effective order conditions and expressed the solutions in terms of the free parameters. 
Then, we have written the corresponding function $ \mathcal{E}_5$ in terms of these parameters and explored
systematically the parameter domain to find its local minima. Finally, we have taken the set of values which provide the smallest value for $ \mathcal{E}_5$. 
In this way we get the coefficient collected in Table \ref{table.4},
with the theoretical efficiencies gathered in Table \ref{tab.effic4}.

\


 \begin{table}[!h]
  \begin{center}
    \begin{tabular}{l|cccc}
      $s$&$E_{\rm{ef}}^{(5)}$&$E_{\rm{ef}}^{(7)}$&$E_{\rm{ef}}^{(9)}$&1--norm\\
      \hline
      3   & 2.2753& 2.5675& 2.6384&4.4048\\
      4   & 1.5470& 1.7567& 1.8304&2.8523\\
      5   & 1.3142& 1.5034& 1.5775&2.3177\\
      6   & 1.2026& 1.3984& 1.4852&2.0417\\
      7   & 1.1389& 1.3220& 1.4207&1.8710\\
      8   & 1.1001& 1.2961& 1.4061&1.7543\\
      9   & 1.0778& 1.2662& 1.3903&1.6672\\      
    \end{tabular}
    \caption{Theoretical efficiencies of $s$-stage kernels $\psi_s^{[4]}$ of effective order 4.}
    \label{tab.effic4}
  \end{center}
 \end{table}

For comparison with the schemes we construct here, the kernel of the processed method BCM$_6^{[4]}$ of Table \ref{table.1.2} 
has the value $E_{\mathrm{ef}}^{(5)} = 1.3432$, whereas the 
corresponding value for the scheme BM$_6^{[4]}$ of Table \ref{table.1} is 1.5829.

Notice that by including more basic maps in the composition it is possible to reduce the value of $E_{\mathrm{ef}}^{(r)}$, and
that errors at
higher order also reduce accordingly. Although several values of the coefficients provide essentially the same (minimum) value for $ \mathcal{E}_5$,
the observed pattern closely follows the rule of thumb formulated by Mclachlan \cite{mclachlan02foh} for kernels of the form
(\ref{eq.2.1.2}) and effective order 4: (i) set the maximum possible number of outer stages in (\ref{eq.2.1.2}) equal to eliminate free parameters; (ii) find
the solution of the effective order conditions with the remaining parameters; (iii) either use the resulting method as it is or take it as the starting point
for further minimize the main error term.

It is shown in \cite{mclachlan02foh} that the minimum effective error is obtained
for $s=19$, which is precisely the same pattern observed here: $E_{\mathrm{ef}}^{(5)}$ achieves its minimum value $\approx 1$ when $s=19$.

Methods with $s=8,9$ in Table \ref{table.4} have been obtained by applying the rule of thumb, but considering two and three free parameters,
respectively, for further optimization.
We notice that, although the value of $E_{\mathrm{ef}}^{(r)}$ diminishes indeed with $s$, it does so more slowly (for instance, for $s=10$ it is
approximately 1.06), so that in practice we restrict ourselves to $s \le 9$.

\begin{table}
  \centering
    \renewcommand\arraystretch{1.1}
    \begin{tabular}{lll}
      \multicolumn{3}{c}{$s=3, \, \psi_{3}^{[4]}$}\\
      \hline 
      $\alpha_1= {\alpha_2} = \frac{1}{6}(2 + 2^{-1/3} + 2^{1/3})$ &\qquad ${\alpha_3}= \frac{1}{2} - 2 \alpha_1$ \\
      \hline 
                                    & & \\
      \multicolumn{3}{c}{$s=4, \, \psi_{4}^{[4]}$}\\
      \hline 
      $\alpha_1= 0.32175$ &\qquad $\alpha_2= -0.46308$ \\
      $\alpha_3= 0.3257797788491147633383644$ & \qquad  $\alpha_4= 0.3155502211508852366616355$ \\
      \hline 
                                    & & \\                                    
      \multicolumn{3}{c}{$s=5, \, \psi_{5}^{[4]}$}\\
      \hline 
      $\alpha_1= \alpha_2 = 0.2014$ &\qquad $\alpha_3= 0.2136$ \\
      $\alpha_4= -0.3294322555468400515907084$ & \qquad  $\alpha_5= 0.2130322555468400515907084$ \\
       \hline 
                                    & & \\                                    
      \multicolumn{3}{c}{$s=6, \, \psi_{6}^{[4]}$}\\
      \hline 
      $\alpha_1= 0.15$ &\qquad $\alpha_2= 0.15$ \\
      $\alpha_3= 0.14353$ & \qquad  $\alpha_4= 0.1592$ \\
      $\alpha_5=-0.26043191662780539180714124$ & \qquad $\alpha_6=0.1577019166278053918071412$ \\
      \hline 
                                    & & \\
      \multicolumn{3}{c}{$s=7, \, \psi_{7}^{[4]}$}\\
      \hline 
      $\alpha_1= 0.1174$ &\qquad $\alpha_2= 0.1158$ \\
      $\alpha_3= 0.1227$ & \qquad  $\alpha_4= 0.112$ \\
      $\alpha_5=0.12685$ & \qquad $\alpha_6=-0.2177553177818524635432753$ \\
      $\alpha_7=0.1230053177818524635432753$ & \\
      \hline 
                                    & & \\
      \multicolumn{3}{c}{$s=8, \, \psi_{8}^{[4]}$}\\
      \hline 
      $\alpha_1= \alpha_2 = \alpha_3 = \alpha_4  = 0.09755$ &\qquad $\alpha_5= 0.09$ \\
      $\alpha_6= 0.1061$ & \qquad  $\alpha_7= -0.1885819261107768579804613$ \\
      $\alpha_8=0.1022819261107768579804613$ & \\
      \hline 
                                    & & \\
     \multicolumn{3}{c}{$s=9, \, \psi_{9}^{[4]}$}\\
      \hline 
      $\alpha_1= \alpha_2 = \cdots = \alpha_7 = 0.082576$ &\qquad  \\
      $\alpha_8=-0.1668033908821750242843527$ & \qquad  $\alpha_9= 0.08877139088217502428435271$ \\
      \hline 
    \end{tabular}
  \caption{\small{Coefficients of time-symmetric kernels of the form (\ref{eq.2.1.1}) of effective order 4 and $s$ stages, $3 \le s \le 9$. \label{table.4}}}
\end{table}

  \subsection{Testing the kernels}
  A simple technique can be used to test the effective order and the relative efficiency of the kernels proposed, before attempting to construct a processor to  form the whole integrator, a task that becomes increasingly difficult with the order.
To proceed, let us consider the linear $d \times d$ matrix system
  \begin{equation} \label{3.2a}  
    \dot{U}=(A_1 + A_2 + A_3)U,\qquad U(0)=I,
  \end{equation}
with exact solution $U_{\rm{ex}}= \e^{t(A_1+ A_2 + A_3)}$, determine the trace of the approximation rendered by the kernel and compare $\tr(\hat{\psi}_h)$ with
$\tr(U_{\rm{ex}})$ at some final time to get an estimate of the error and efficiency that this particular kernel can achieve by processing. 
We take $A_j$, $j=1,2,3$, as $50 \times 50$ matrices whose entries are chosen
randomly from a normal distribution. In addition, to illustrate
the role that the basic scheme may play in the overall performance of the method, we take two different choices for $\chi_h$, namely a composition of the exact flow
of each sub-part,
\begin{equation} \label{exc}
  \chi_h = \e^{h A_3} \, \e^{h A_1} \, \e^{h A_1}
\end{equation}
and the first-order approximation
\begin{equation} \label{foa}
  \chi_h = (I + h A_3)(I + h A_2)(I + h A_1),
\end{equation}
in which case $\chi_h^* = (I - h A_1)^{-1} (I - h A_2)^{-1} (I - h A_3)^{-1}$. 
Figure \ref{fig:tracesa} corresponds to one particular instance of the first case, and Figure \ref{fig:tracesb} corresponds to the second. Here kernels of Table
\ref{table.4} are tested, together with BM$_6^{\blau{[4]}}$ and the kernel of method BCM$_6^{[4]}$ (Table \ref{table.1.2}). Whereas there is
a good deal of variability in the efficiency exhibited by the different schemes with the choice \eqref{exc} depending on the particular matrices $A_j$,
this is not the case of \eqref{foa}: in all examples we have tested, the efficiency closely follows the pattern shown in Table \ref{tab.effic4}.

 \begin{figure}
  \begin{center}
    \subfloat[]{
      \label{fig:tracesa}
      \includegraphics[width=8cm]{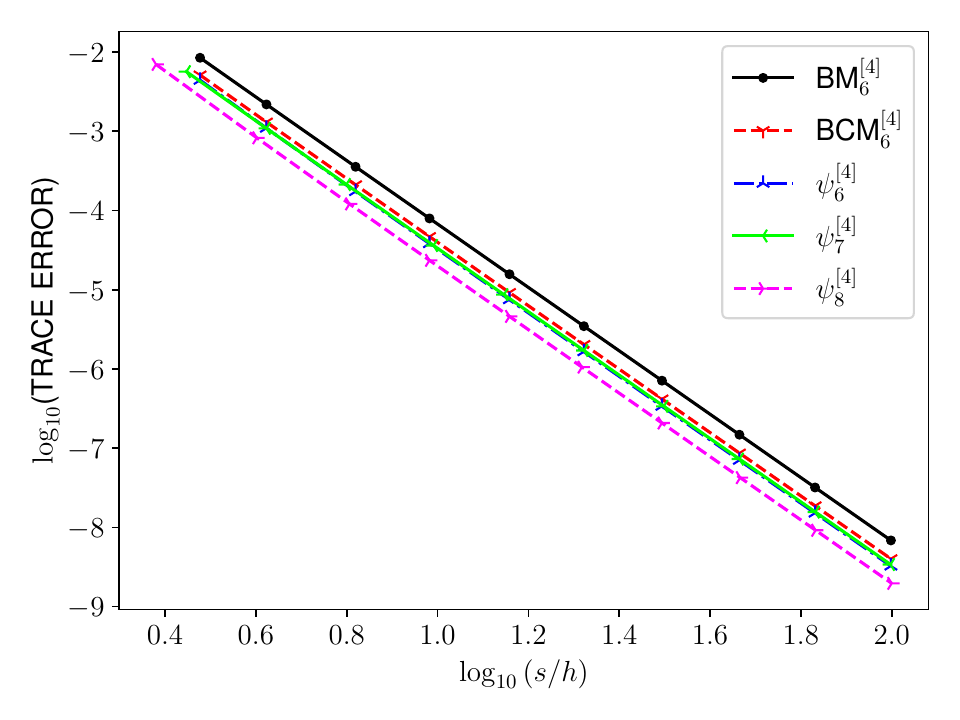}}
    \subfloat[]{
      \label{fig:tracesb}
      \includegraphics[width=8cm]{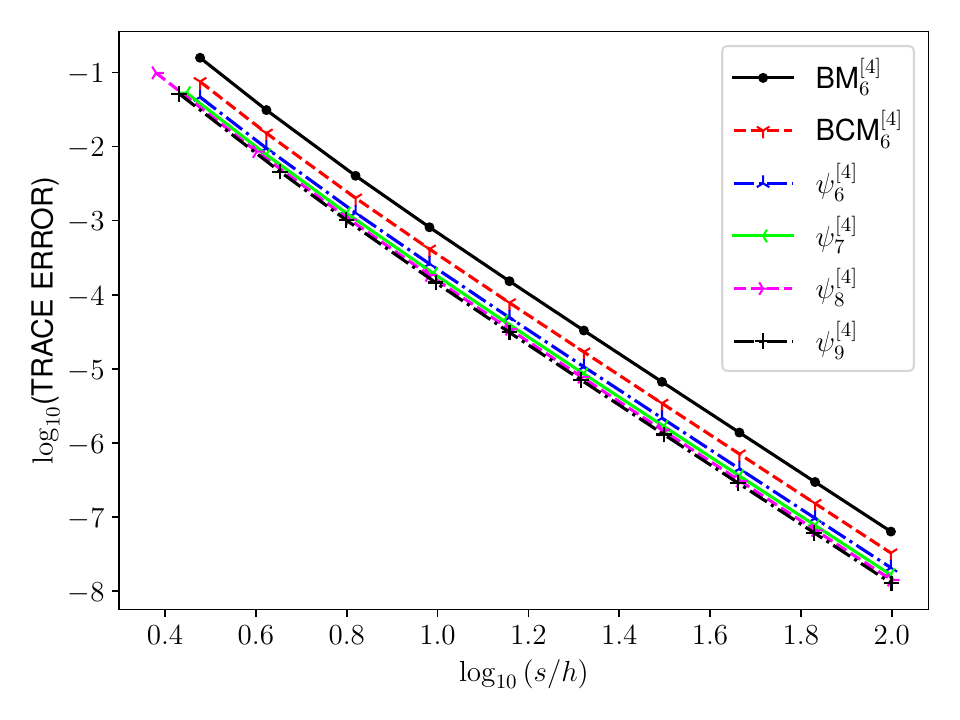}}
    \caption{Error in the trace obtained by kernels of effective order 4 collected in Table \ref{table.4}, together with method BM$_6^{[4]}$ and the kernel of
      BCM$_6^{[4]}$ when they are applied to the linear system \eqref{3.2a} at final time $t_f=10$. (a): basic scheme $\chi_h$ taken as \eqref{exc}; (b): $\chi_h$
      given by \eqref{foa}.}
    \label{fig:traces}
  \end{center}
\end{figure}

\subsection{Construction of the processor}

To have a complete integrator we must obtain a processor $\pi_h$ in (\ref{proces.1}) once a particular kernel has been chosen. This in principle
would require determining the exact flow of the infinite series (\ref{p.series}). For an easier implementarion, however, 
it is more convenient to construct an approximation to $\pi_h$
as a composition of the same form as the kernel. In our case it would be enough to fix conditions on $p_{i,j}$ in 
(\ref{p.series}) up to order 3, but the overall error is reduced if in addition all the conditions are satisfied up to order 5, as previously explained. 
In other words, we fix
\[
\begin{aligned}
   & p_{1,1} = p_{3,1} = p_{3,2} = 0, \qquad p_{2,1} = {k_{3,2}}, \\
   & p_{4,1} = {k_{5,2}}, \qquad p_{4,2} = {k_{5,3}}, \qquad p_{4,3} = {k_{5,4}}
\end{aligned}
\]   
for a particular kernel, and then construct an approximation of the form
\begin{equation}  \label{proc.compo}
  \pi_h \approx  \pi^{(s,4)} \equiv
\psiba_{\beta_7h}\circ  \psib_{\beta_6 h} \circ \chi_{\beta_5 h}^* \circ \chi_{\beta_4 h} \circ \chi_{\beta_3 h}^*
\circ \chi_{\beta_2 h} \circ \chi_{\beta_1 h}^*,
\end{equation}
assuming that the corresponding equations have real solutions. Then, we approximate the inverse map $\pi_h^{-1}$ by
\[
  \pi_h^{-1} \approx  
\psib_{-\beta_1h}\circ  \psiba_{-\beta_2 h} \circ \chi_{-\beta_3 h} \circ \chi_{-\beta_4 h}^* \circ \chi_{-\beta_5 h}
\circ \chi_{-\beta_6 h}^* \circ \chi_{-\beta_7 h}.
\] 
As explained in \cite{blanes21sps,blanes23esm}, one can also replace $\pi_h^{-1}$ by the adjoint $\pi_h^*$ and obtain an approximation up to the same order.
 In that case, the integrator after $N$ steps reads
\[
  \pi_h \circ \psi_h^N \circ \pi_h^*
\]
and is also time-symmetric if $\psi_h$ is so, although it cannot be obtained as a $N$-fold composition of a one-step map \cite{blanes21sps}.
In consequence, we take
\[
  \pi_h^* \approx  
\psib_{\beta_1h}\circ  \psiba_{\beta_2 h} \circ \chi_{\beta_3 h} \circ \chi_{\beta_4 h}^* \circ \chi_{\beta_5 h}
\circ \chi_{\beta_6 h}^* \circ \chi_{\beta_7 h}.
\] 

By applying this strategy, we have determined a processor for each kernel of Table \ref{table.4} (except for $\psi_3$, given its poor efficiency). 
The overall methods read as
 \begin{equation}\label{eq.psihat4}
    \hat{\psi}^{[s,4]} = \pi^{(s,4)} \circ \psi_s^{[4]} \circ \left(\pi^{(s,4)}\right)^{*},
  \end{equation}
where $s$ denotes the number of stages of the kernel.

\section{Processed methods of order 6}
\label{sect4}
A similar procedure can be carried out to construct processed composition methods of order 6. In this case one has 5 effective order
conditions, so that the kernel involves at least $s=5$ stages. Methods with $s \le 10$ have been obtained in \cite{blanes06cmf} by applying the previous rule
of thumb. They are more efficient than the standard scheme BM$_{10}^{[6]}$ of Table \ref{table.1} on a number of examples. Proceeding
analogously as in the case of order 4, we build
 new kernels by minimizing the objective function $E_{\mathrm{ef}}^{(7)}$. This is done by exploring
the space of
parameters and identifying local minima of $\mathcal{E}_7$  for $5 \le s \le 11$. The coefficients of the most promising schemes we have found 
are collected in   
Table \ref{table.5}, leading to the following efficiencies shown in Table \ref{effic_6}.

\


 \begin{table}[!h]
  \begin{center}
    \begin{tabular}{l|ccc}
      $s$&$E_{\rm{ef}}^{(7)}$&$E_{\rm{ef}}^{(9)}$&1--norm\\
      \hline
      5& 5.1053& 6.2573& 9.6024\\
      6& 3.1347& 3.6818& 5.7329\\
      7& 2.5170& 2.9862& 4.3759\\
      8& 2.2193& 2.6388& 3.6553\\
      9& 2.0686& 2.4786& 3.2417\\
     10& 1.9488& 2.3345& 2.9099\\
     11& 1.8718& 2.2560& 2.6935\\
    \end{tabular}
    \caption{Theoretical efficiencies of $s$-stage kernels $\psi_s^{[6]}$ of effective order 6.}
    \label{effic_6}
  \end{center}
   \end{table}

 We should stress that the kernel with $s=5$ stages also corresponds to a composition (\ref{eq.2.1.2}).
For comparison, the standard 6th-order method BM$_{10}^{[6]}$ of Table \ref{table.1} has
$E_{\mathrm{ef}}^{(7)} = 3.5855$,
whereas the kernel BCM$_9^{[6]}$ has
$E_{\mathrm{ef}}^{(7)} = 2.2144$. It is worth remarking that the most efficient schemes we have found with $s=10$ and $s=11$ correspond precisely to a particular case of the rule of thumb. Moreover, these schemes 
correspond precisely to compositions of the form \eqref{eq.2.1.2}. A careful exploration of the region of free parameters
has not provided solutions with a better efficiency.


Concerning the processor, one has to solve 23 equations to achieve order 6, so we approximate $\pi_h$ by a composition of the form (\ref{proc.compo}), 
$\pi^{(s,6)} = \chi_{\beta_{23}h}^* \circ \chi_{\beta_{22} h} \circ  \cdots \circ \chi_{\beta_{1}h}^*$. Among all the real solutions obtained, we take the one with the minimum 1-norm  of the vector $(\beta_1, \ldots, \beta_{23})$.
The overall method is now denoted as 
  \begin{equation}\label{eq.psihat6}
    \hat{\psi}^{[s,6]} = \pi^{(s,6)} \circ \psi_s^{[6]} \circ \left(\pi^{(s,6)}\right)^{*}
  \end{equation}
  and methods with $5 \le s \le 11$ have been constructed.

We have also explored kernels of order 8, but all the schemes we have been able to construct correspond indeed to compositions of a time-symmetric second
order basic method.

\begin{table}
  \centering
    \renewcommand\arraystretch{1.1}
    \begin{tabular}{lll}
      \multicolumn{3}{c}{$s=5$, $\psi_{5}^{[6]}$}\\
      \hline 
      $\alpha_1= 1.1983882307745148$ &\qquad $\alpha_2= -1.0753056449710827$ \\
      $\alpha_3= -1.0753056449710827$ & \qquad $\alpha_4 = 0.7261115295838254$ \\
      $\alpha_5= 0.7261115295838252$ &   \\
      \hline 
                                    & & \\      
      \multicolumn{3}{c}{$s=6$, $\psi_{6}^{[6]}$}\\
      \hline 
      $\alpha_1= 0.3579656411737745366807624$ &\qquad $\alpha_2= 0.3041155195721355055671815$ \\
      $\alpha_3= 0.3544845132692152058236740$ & \qquad $\alpha_4 = -0.5776359154029903584192191$ \\
      $\alpha_5= -0.6055964252788015563190656$ & \qquad $\alpha_6 = \frac{2}{3}$ \\
      \hline 
                                    & & \\
      \multicolumn{3}{c}{$s=7$, $\psi_{7}^{[6]}$}\\
      \hline 
      $\alpha_1= 0.2$ &\qquad $\alpha_2=0.2102$ \\
      $\alpha_3= 0.2076682089468184822550893$ & \qquad  $\alpha_4= 0.2483663566422618080555971$ \\
      $\alpha_5= -0.4108957823061925870283283$ & \qquad $\alpha_6 = -0.4330744093869198010728686$ \\
      $\alpha_7= 0.477735626104032097790510$ & \\
      \hline 
                                    & & \\                                    
      \multicolumn{3}{c}{$s=8$, $\psi_{8}^{[6]}$}\\
      \hline 
      $\alpha_1= 0.1535$ &\qquad $\alpha_2= 0.146$ \\
      $\alpha_3= 0.1535$ & \qquad  $\alpha_4 = 0.1564865138360775523331602$ \\
      $\alpha_5= 0.1777546764340215024573463$ & \qquad $\alpha_6 = -0.3260392072026447259933467$ \\
      $\alpha_7= -0.3377852074639320941003920$ & \qquad $\alpha_8 = 0.376583224396477765303232$ \\
       \hline 
                                    & & \\                                    
      \multicolumn{3}{c}{$s=9$, $\psi_{9}^{[6]}$}\\
      \hline 
      $\alpha_1= 0.1145$ &\qquad $\alpha_2= 0.116$ \\
      $\alpha_3= 0.117$ & \qquad  $\alpha_4= 0.1115$ \\
      $\alpha_5=0.1319890385474291590082355$ & \qquad $\alpha_6=0.1512264299418584439400894$ \\
      $\alpha_7= -0.2763628586973694504837428$ & \qquad $\alpha_8 = -0.2840658003186325771164088$ \\
      $\alpha_9 = 0.318213190526714424651826$ & \\
       \hline 
                                    & & \\                                    
      \multicolumn{3}{c}{$s=10$, $\psi_{10}^{[6]}$}\\
      \hline 
      $\alpha_1= \cdots=\alpha_7= 0.100838384835000970361478569216$ &\\
      $\alpha_8=\alpha_9=-0.238737866770265639777320334936 $ &\\
      $\alpha_{10}= 0.271607039695524487024290685363$&\\
      \hline
                                     & & \\      
            \multicolumn{3}{c}{$s=11$, $\psi_{11}^{[6]}$}\\
      \hline 
      $\alpha_1=\cdots=\alpha_8= 0.0852884432504611078508$ &\\
      $\alpha_9=\alpha_{10}=  -0.2116830704463290239945$ &\\
      $\alpha_{11}= 0.241058594888969185183038787789$&\\
      \hline 
     \end{tabular}
  \caption{\small{Coefficients of time-symmetric kernels of the form (\ref{eq.2.1.1}) of effective order 6 with $s$ stages, $5 \le s \le 11$. \label{table.5}}}
\end{table}

\section{Numerical experiments}

Next we test some of the previous processed methods for the numerical integration of the equations describing the
motion of a particle in an electromagnetic field, and the motion
around a Reissner--Nodstr\"om black hole. In the first case, the system is split into three parts, whereas in the second one has to consider five parts.
In this respect, as noted in \cite{mclachlan22tsi}, choosing the right sequence of basic maps when forming $\chi_h$ is of the utmost importance. In fact,
the numerical experiments carried out in \cite{mclachlan22tsi} show that the overall efficiency of the scheme may vary significantly with the ordering. For
this reason, it is recommended to test different sequences to identify the most efficient for a given problem.

Several codes implementing the previous schemes for the two problems at hand are available at

\begin{center}
  doi.org/10.5281/zenodo.8375196.
\end{center}

Specifically, the interested reader can find there the coefficients $\beta_j$ of the processor for methods of order 4 and 6, eqs. \eqref{eq.psihat4} and  
\eqref{eq.psihat6}, together with two new processors for the 4th- and 6th-order proposed in \cite{blanes06cmf} and collected in Table \ref{table.1.2} and
the codes generating the figures included in this section.

\subsection{Motion of a charged particle under Lorentz force}\label{sec:lor}

The evolution of a  particle of mass $m$ and charge $q$ in a external electromagnetic field (in the non-relativistic limit) is modeled by
the equation
\begin{equation}   \label{lorentz.1}
   m \, \ddot{\mathbf{x}} = q \, (\mathbf{E} + \dot{\mathbf{x}} \times \mathbf{B}),
\end{equation}   
which can also be written as a first-order system:
\begin{equation}   \label{lorentz.2}
 \begin{aligned}
   &  \dot{\mathbf{x}} = \mathbf{v} \\
   &  \dot{\mathbf{v}} = \frac{q}{m}  \mathbf{E} + \omega \, \mathbf{b} \times \mathbf{v}.
  \end{aligned}
\end{equation}
Here $\omega = -q B/m$ is the local cyclotron frequency, $B = \|\mathbf{B}\|$ and $\mathbf{b} = \mathbf{B}/B$ is the unit vector in the direction of
the magnetic field. In the following, we assume that both $ \mathbf{E}$ and $ \mathbf{B}$ depend only on the position variables $\mathbf{x}$. 

System (\ref{lorentz.2}) can be split into three parts in such a way that each subpart is explicitly solvable, and preserve volume in phase space
$(\mathbf{x}, \mathbf{v})$  \cite{he15vpa,he16hov}. Specifically, if we define
$y=(\mathbf{x}, \mathbf{v})^\top$, then \eqref{lorentz.2} can be expressed as
\begin{eqnarray}   \label{lorentz.3}
  \dot{y} = \frac{d}{dt} \left(  \begin{array}{c}
  				\mathbf{x}  \\
				\mathbf{v} 
			  \end{array}  \right)  & = &  
		 \left(  \begin{array}{c}
  				\mathbf{v}  \\
				0 
			  \end{array}  \right)	  + 
		 \left(  \begin{array}{c}
  				0  \\
				\frac{q}{m} \mathbf{E}(\mathbf{x}) 
			  \end{array}  \right)	  	+ 
		 \left(  \begin{array}{c}
  				0  \\
				\omega(\mathbf{x})  \mathbf{b}(\mathbf{x}) \times \mathbf{v}
			  \end{array}  \right),	  \nonumber 	\\
                          & = & f^{[A]}(y) +  f^{[B]}(y) + f^{[C]}(y),
\end{eqnarray}
with exact solutions given by
\begin{equation}   \label{lorentz.4}
  \varphi_t^{[A]}: \left\{  \begin{array}{l}
   					\mathbf{x}(t) = \mathbf{x}_0 + t \, \mathbf{v}_0 \\
					\mathbf{v}(t) = \mathbf{v}_0
				\end{array}  \right. \qquad 	
   \varphi_t^{[B]}: \left\{  \begin{array}{l}
   					\mathbf{x}(t) = \mathbf{x}_0  \\
					\mathbf{v}(t) = \mathbf{v}_0 + t \, \frac{q}{m} \, \mathbf{E}(\mathbf{x}_0) 
				\end{array}  \right.  	\qquad
   \varphi_t^{[C]}: \left\{  \begin{array}{l}
   					\mathbf{x}(t) = \mathbf{x}_0  \\
					\mathbf{v}(t) = \e^{t \omega_0 \hat{\mathbf{b}}_0} \,  \mathbf{v}_0
				\end{array}  \right.  
\end{equation}
in terms of $\omega_0 \equiv \omega(\mathbf{x}_0)$ and $\hat{\mathbf{b}}_0 \equiv \hat{\mathbf{b}}(\mathbf{x}_0)$, with 
\[
   \hat{\mathbf{b}}(\mathbf{x}) = \left(  \begin{array}{ccc}
   				0  &  -b_3(\mathbf{x})  &  b_2(\mathbf{x}) \\
			  b_3(\mathbf{x})   &  0  &  -b_1(\mathbf{x}) 	\\
			  -b_2(\mathbf{x})  &  b_1(\mathbf{x}) &  0
			      \end{array}  \right).
\]
In practice, as in \cite{he15vpa}, we use the expression
\[
  \exp(t \omega_0 \hat{\mathbf{b}}_0)  \mathbf{v}_0 = \mathbf{v}_0+\sin(t\omega_0)\hat{\mathbf{b}}_0\mathbf{v}_0
  +(1-\cos(t\omega_0))\hat{\mathbf{b}}_0^2\mathbf{v}_0
\]
for computing $ \varphi_t^{[C]}$, and consider
a static, non-uniform electromagnetic field
\begin{equation}   \label{loretnz.5}
   \mathbf{E} = - \nabla V = \frac{\alpha}{r^3} (x \, \mathbf{e}_x + y \ \mathbf{e}_y), \qquad   \mathbf{B} = \nabla \times \mathbf{A} = r \, \mathbf{e}_z,
\end{equation}
derived from the potentials
\[
   V = \frac{\alpha}{r}, \qquad\qquad \mathbf{A} = \frac{r^2}{3} \, \mathbf{e}_{\theta},
\]
respectively, in cylindrical coordinates $(r, \theta, z)$. The parameter $\alpha$ is used to parametrize the scalar potential. Then, both the angular momentum
and energy
\[
   L = r^2 \dot{\theta} + \frac{r^3}{3}, \qquad \qquad 
    H = \frac{1}{2} \|\mathbf{v}\|^2 + \frac{\alpha}{r}
\]
are invariants of the problem. 

For the simulations we take $q=-1$, $m=1$, initial position $\mathbf{x}_0 = (0, -1, 0)^T$ and initial velocity $\mathbf{v}_0 = (0.1, 0.01, 0)$ \cite{he15vpa}. In our
first experiment we integrate until the final time $t_f = 20$ with the 4th-order processed methods $\hat{\psi}^{[6,4]}$ and $\hat{\psi}^{[8,4]}$ 
(c.f. eqs. \eqref{eq.psihat4} and  \eqref{eq.psihat6})
and values
of $\alpha$ in the interval $(0,1/10)$. As basic scheme we choose the Lie--Trotter splitting
\begin{equation} \label{LT1}
 \chi_h =  \varphi_h^{[A]} \circ \varphi_h^{[B]} \circ \varphi_h^{[C]}.
\end{equation} 
We determine the error
in phase space at the final time with each integrator by taking as reference solution the output generated by the standard routine DOP853 \cite{hairer93sod}. The results are depicted in Figure \ref{fig:em1a}, where the errors committed 
by the standard method BM$_6^{[4]}$ and the processed scheme BCM$_6^{[4]}$ are also shown. The step size $h$ is chosen in such a way that all integrators
require the same computational cost. 

Figure \ref{fig:em1b} corresponds to an efficiency diagram obtained by methods BM$_6^{[4]}$, BCM$_6^{[4]}$ and the new processed schemes $\hat{\psi}^{[6,4]}$
and $\hat{\psi}^{[8,4]}$
when $\alpha = 0.07$ at the final time $t_f = 200$. Notice that the processed schemes provide more accurate results with the same computational effort.
For comparison, we have also depicted the result achieved by the triple jump composition \eqref{yoshida}, SS$_3^{[4]}$.

Figures \ref{fig:em2a} and \ref{fig:em2b} show the results achieved, for the same problem, initial conditions and final integration time, by the following
6th-order schemes: BM$_{10}^{[6]}$ (Table \ref{table.1}), BCM$_9^{[6]}$ (Table \ref{table.1.2} with the new processor), $\hat{\psi}^{[8,6]}$ and $\hat{\psi}^{[10,6]}$. Here again
the new processed scheme is more efficient when high accuracy is desired. We also include for comparison the result achieved by the most efficient 7-stage
6th-order symmetric composition of 2nd-order schemes  proposed in \cite{yoshida90coh}, SS$_7^{[6]}$.

Finally, Figure \ref{fig:em_ord} shows the non-trivial effects on the overall
error of the different orderings in the basic scheme $\chi_h$  for several values of $\alpha$.

\begin{figure}
  \begin{center}
    \subfloat[]{
      \label{fig:em1a}
      \includegraphics[width=8cm]{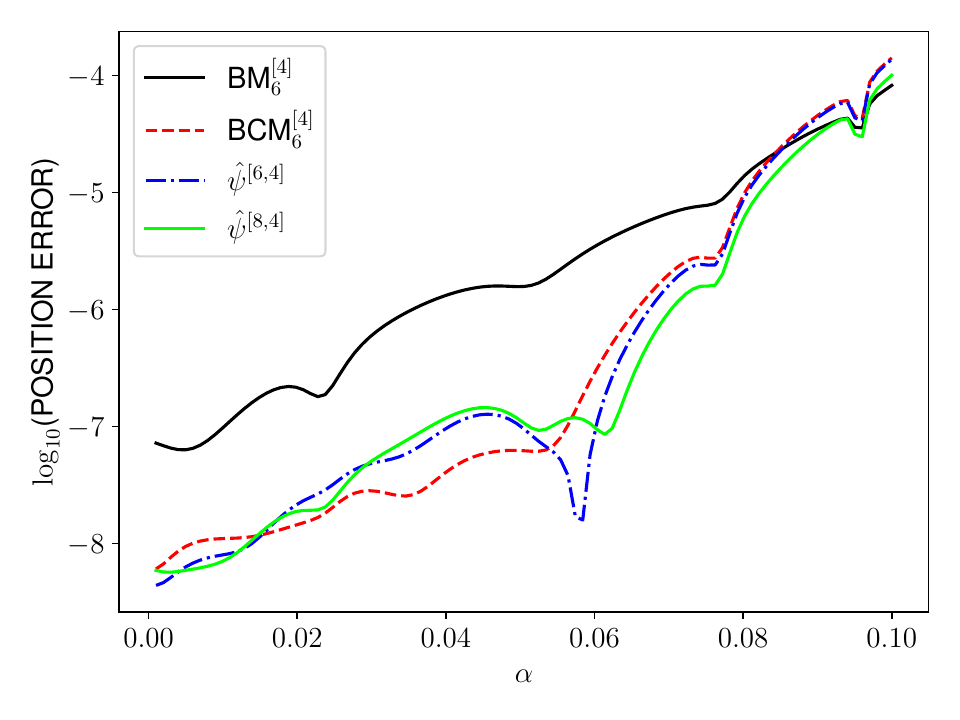}}
    \subfloat[]{
      \label{fig:em1b}
      \includegraphics[width=8cm]{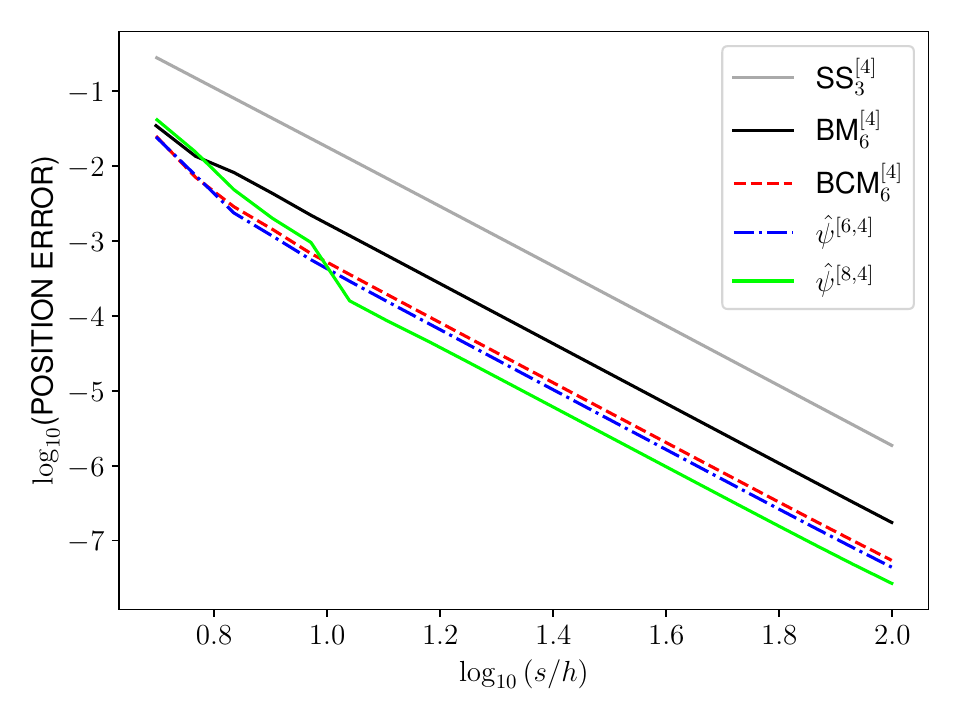}}
    \caption{\textbf{Motion of a charged particle under Lorentz force}. (a) Maximum position error for different values of the parameter $\alpha$ and for various 4th-order methods at the final time $t_f=200$, evaluated with the same computational cost $s/h=40$.
    (b) Efficiency diagram for different 4th-order methods with $\alpha=0.07$ and $t_f=200$.}
    \label{fig:em1}
  \end{center}
\end{figure}

\begin{figure}
  \begin{center}
    \subfloat[]{
      \label{fig:em2a}
      \includegraphics[width=8cm]{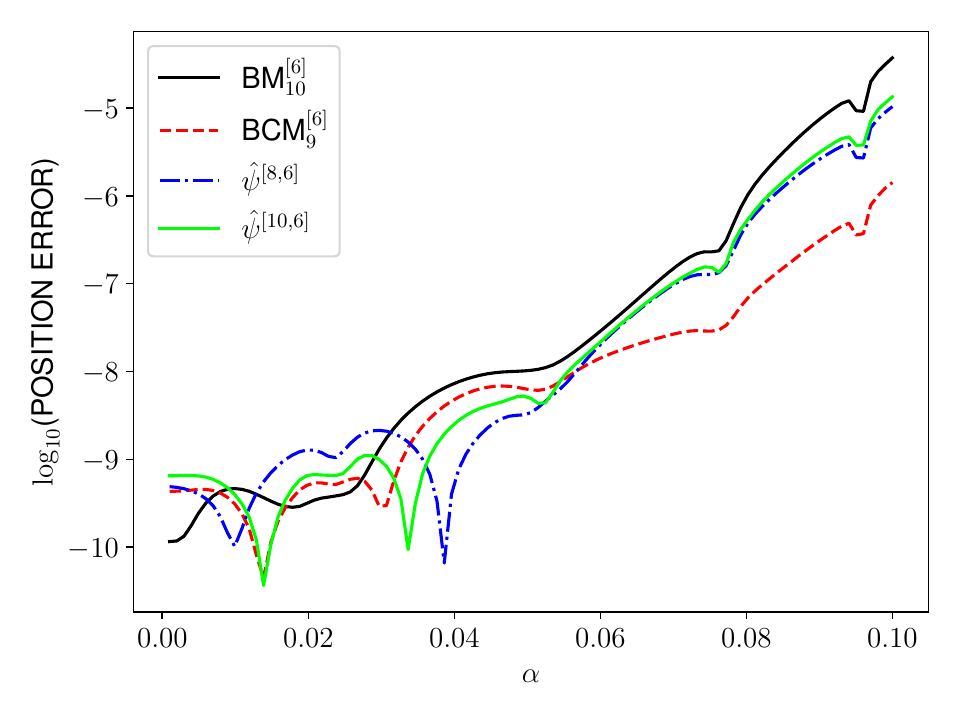}}
    \subfloat[]{
      \label{fig:em2b}
      \includegraphics[width=8cm]{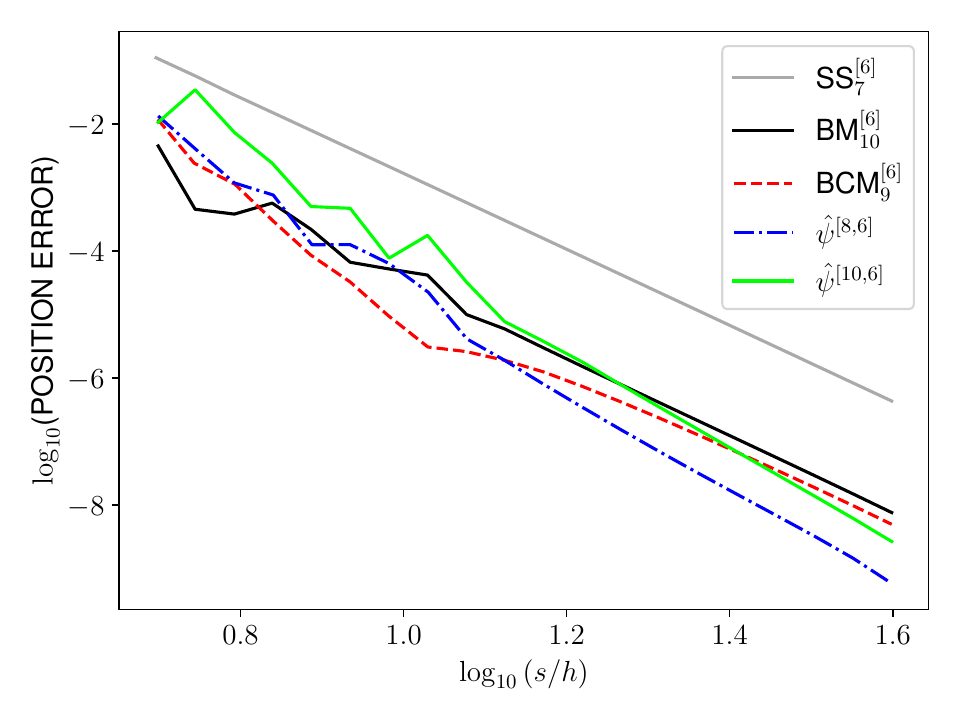}}
    \caption{\textbf{Motion of a charged particle under Lorentz force}. (a) Maximum position error for different values of the parameter $\alpha$ and for various 6th-order methods at the final time $t_f=200$, evaluated with the same computational cost $s/h=40$.
    (b) Efficiency diagram for different 6th-order methods with $\alpha=0.04$ and $t_f=200$.}
    \label{fig:em2}
  \end{center}
\end{figure}

\begin{figure}
  \begin{center}
    \includegraphics[width=10cm]{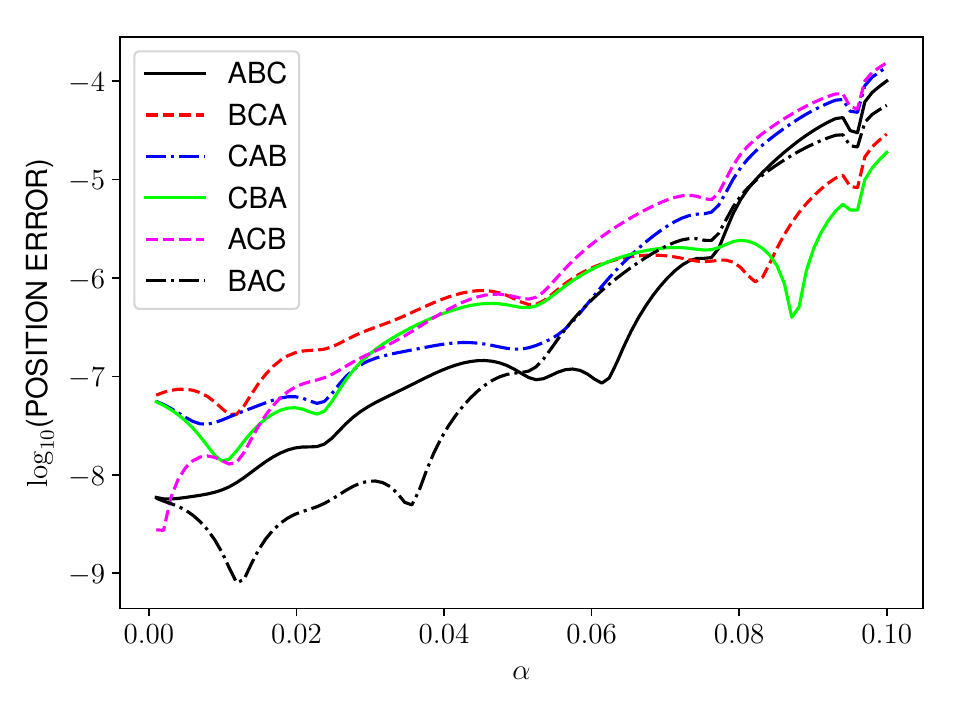}
    \caption{\textbf{Motion of a charged particle under Lorentz force}.
      Maximum position error for the method $\hat{\psi}^{[8,4]}$ for different values of the parameter $\alpha$ with all possible orderings of the first-order method $\chi_h$. The simulation was conducted with $s/h=40$ and $t_f=200$.}
    \label{fig:em_ord}
  \end{center}
\end{figure}

  \subsection{Particle around a Reissner--Nordström black hole}
  A Schwarzschild black hole with charge $Q$ is known as a Reissner--Nordstr\"om black hole. 
  The motion of a test particle around this black hole is described by the Hamiltonian \cite{wang21coe2}:
  \begin{equation}\label{rnbh.ham}
    H =\frac{-1}{2\left(1-\frac{2}{r}+ \frac{Q^2}{r^2}\right)}E^2+\frac{1}{2}\left(1-\frac{2}{r} + \frac{Q^2}{r^2}\right)p_r^2+\frac{1}{2}\frac{p_\theta^2}{r^2}+ \frac{L^2}{2r^2\sin^2\theta}.
  \end{equation}
  Here $r$ and $\theta$  correspond to the radial and angular coordinates of the particle, $p_r$ and $p_\theta$ are their conjugate momenta, and $E$ and $L$ are the energy and angular momentum of the particle, respectively.

  This Hamiltonian can be separated into five explicitly integrable parts, namely \cite{wang21coe2}
  \begin{equation*}
    H=H_A+H_B+H_C+H_D+H_E,
  \end{equation*}
  with
  \[
  \begin{aligned}
    & H_A=\frac{-1}{2\left(1-\frac{2}{r}+ \frac{Q^2}{r^2}\right)}E^2+ \frac{L^2}{2r^2\sin^2\theta},\qquad\qquad H_B=\frac{1}{2}p_r^2,\\
    & H_C=-\frac{1}{r}p_r^2, \qquad\qquad 
     H_D=\frac{p_\theta^2}{2r^2},\qquad\qquad
    H_E=\frac{Q^2p_r^2}{2r^2},
  \end{aligned}
\]  
whose flows read explicitly
\begin{equation} \label{sbh.flows}
\begin{aligned}   
    &  \varphi_\tau^{[A]}: \left\{
       \begin{array}{l}
         r(\tau) = r_0 \\
         p_r(\tau) = {p_r}_0 + \tau\left(\frac{L^2}{r_0^3\sin^2\theta_0}-\frac{E^2}{(r_0-2)^2}\right)\\
         \theta(\tau) = \theta_0\\
         p_\theta(\tau) = {p_\theta}_0 + \tau\frac{L^2\cos\theta_0}{r_0^2\sin^3\theta_0}
       \end{array}  \right. \qquad\qquad	
  \varphi_\tau^{[B]}: \left\{
  \begin{array}{l}
         r(\tau) = r_0 +\tau {p_r}_0\\
         p_r(\tau) = {p_r}_0 \\
         \theta(\tau) = \theta_0\\
         p_\theta(\tau) = {p_\theta}_0    			
  \end{array}  \right. \\[0.2em]
   &  \varphi_\tau^{[C]}: \left\{
       \begin{array}{l}
         r(\tau) = \left(\frac{(r_0^2-3\tau{p_r}_0)^2}{r_0}\right)^{1/3} \\
         p_r(\tau) = {p_r}_0\left(\frac{(r_0^2-3\tau{p_r}_0)}{r_0^2}\right)^{1/3} \\
         \theta(\tau) = \theta_0\\
         p_\theta(\tau) = {p_\theta}_0  
       \end{array}  \right.\;\;
  \varphi_\tau^{[D]}: \left\{
  \begin{array}{l}
    r(\tau) = r_0 \\
    p_r(\tau) = {p_r}_0+\tau\frac{{p_\theta}^2_0 }{r_0^3} \\
    \theta(\tau) = \theta_0+\tau\frac{{p_\theta}_0 }{r_0^2}\\
    p_\theta(\tau) = {p_\theta}_0   			
  \end{array}  \right. \quad 
  \varphi_\tau^{[E]}: \left\{
       \begin{array}{l}
         r(\tau) =  r_0 + \tau {p_r}_0 \frac{Q^2}{r^2}\\[0.2em]
         p_r(\tau) =  {p_r}_0 + \tau \frac{Q^2{p_r}_0^2}{r^3}\\
         \theta(\tau) = \theta_0\\
         p_\theta(\tau) = {p_\theta}_0  
       \end{array}  \right.	  		
  \end{aligned}
\end{equation}

\begin{figure}
  \begin{center}
    \subfloat[]{
      \label{fig:rn1a}
      \includegraphics[width=8cm]{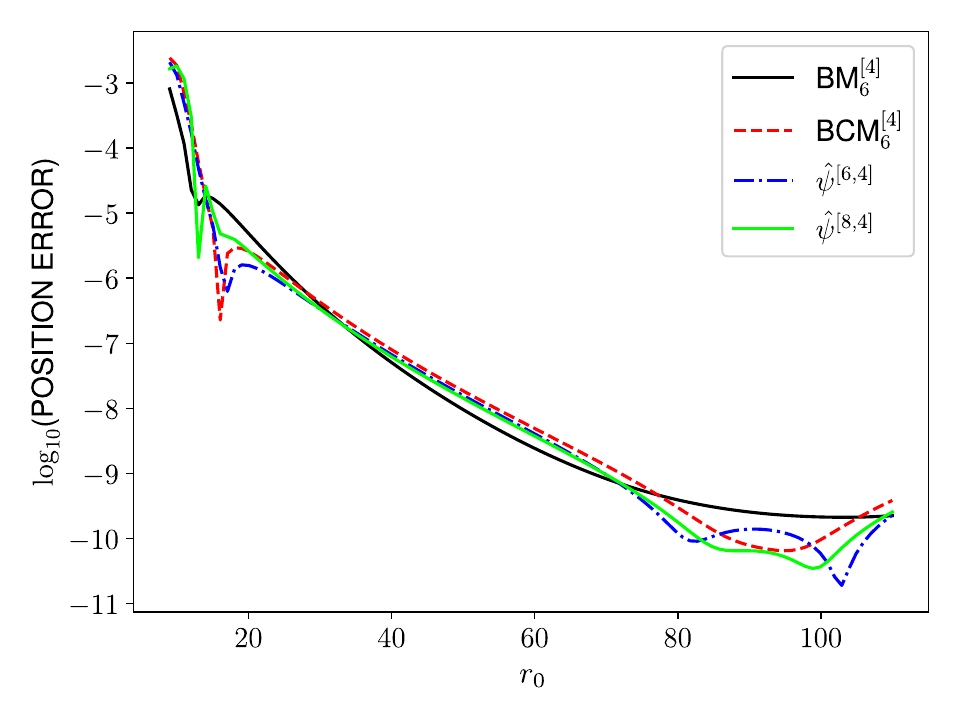}}
    \subfloat[]{
      \label{fig:rn1b}
      \includegraphics[width=8cm]{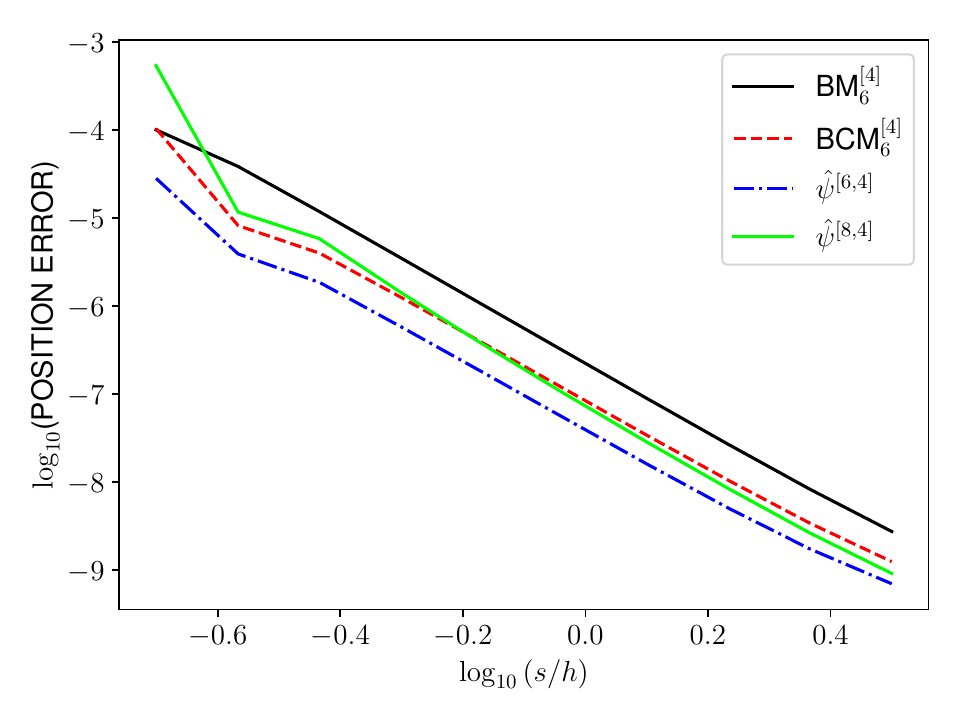}}
    \caption{\textbf{Particle around Reissner--Nodström black hole}. (a) Maximum position error for different values of the coordinate $r_0$ and for various 4th-order methods  at the final time $t_f=10^4$, evaluated with the same computational cost $s/h=0.4$
    (b) Efficiency diagram for different 4th-order methods with $r_0=18$ and the same $t_f$.}
    \label{fig:rn1}
  \end{center}
\end{figure}

We take $E=0.995$, $L=4.6$, $Q=0.3$, initial angle $\theta_0=\pi/2$,  and integrate until the final time $t_f=10^4$ by taking as
basic scheme $\chi_h=\varphi_h^{[A]}\circ\varphi_h^{[B]}\circ\varphi_h^{[C]}\circ\varphi_h^{[D]}\circ\varphi_h^{[E]}$. In Figure \ref{fig:rn1a} we depict
the error in phase space by varying $r_0$ in the interval $r_0 \in[9,110]$ for methods BM$_6^{[4]}$, BCM$_6^{[4]}$ and the new
  processed methods $\hat{\psi}^{[6,4]}$ and $\hat{\psi}^{[8,4]}$. The step size is taken so that all the methods require the same computational effort.
The right panel \ref{fig:rn1b} shows the corresponding efficiency diagram, obtained with the same initial condition and $r_0=18$. 

\begin{figure}
  \begin{center}
    \subfloat[]{
      \label{fig:subfig-1}
      \includegraphics[width=8cm]{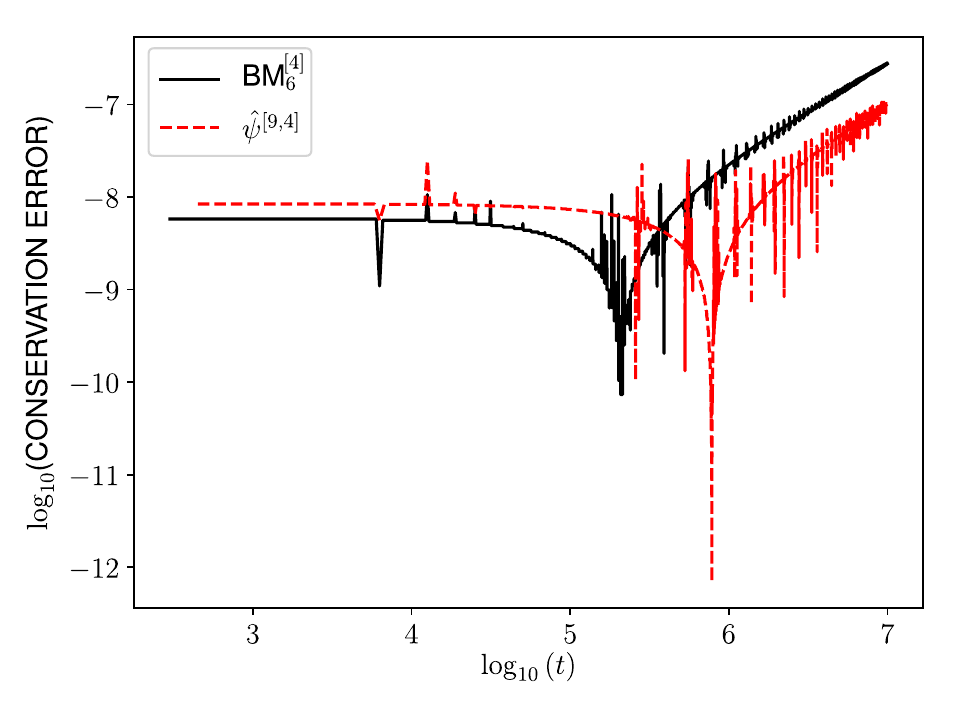}}
    \subfloat[]{
      \label{fig:subfig-2}
      \includegraphics[width=8cm]{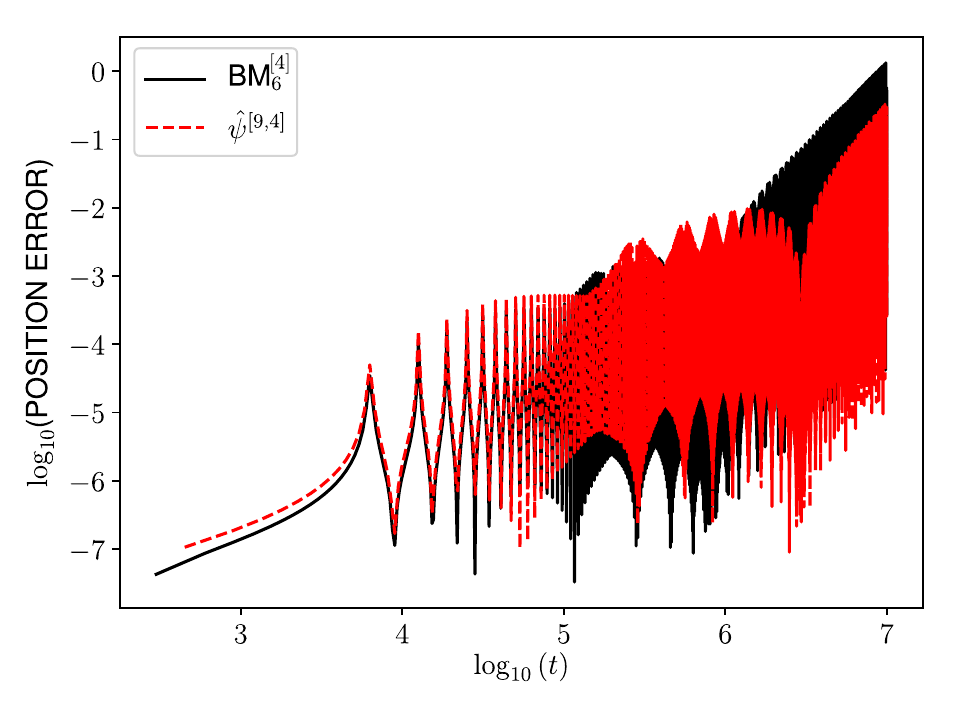}}
    \caption{\textbf{Particle around Reissner--Nodström black hole}. The evolution of the energy error (left) and position error (right) for BM$_6^{[4]}$ and 
    $\hat{\psi}^{[9,4]}$. Both simulations are conducted with the same computational cost $s/h=1$ and a step size $h<10$.}
    \label{fig:rnZ}
  \end{center}
\end{figure}
In addition to the improvement in the efficiency of processing methods with respecto to standard compositions in this case, there is an important aspect to highlight. In \cite{wang21coe1}, it is claimed that in these types of problems, roundoff errors grow and eventually lead to a drift in the conservation of the Hamiltonian when $h>10$. 

Our tests show that this phenomenon, although still present for processed methods, is  delayed with respect to standard compositions, 
as illustrated in Figure \ref{fig:rnZ}.

  \subsection{Motion of a charged particle under Lorentz force II}

The effect of the ordering in the basic method $\chi_h$ is clearly visible e.g. in Figure \ref{fig:em_ord}. We observe that, for some values of the parameter
$\alpha$, the error is up to 100 times smaller for the same scheme, and this may somehow conceal the potential advantages of one particular method
with respect to others. 

We should take into account, however, that the methods considered in this work are all based in compositions of an arbitrary first-order scheme $\chi_h$
and its adjoint $\chi_h^*$, and not just on the Lie--Trotter splitting \eqref{LT1}. It makes sense, then, to analyze the relative performance of the different methods
on a given problem for another choice of $\chi_h$. To do that, we consider again the problem (\ref{lorentz.2}) and take the explicit Euler method as $\chi_h$,
so that $\chi_h^*$ is given by the implicit Euler method. Figure \ref{fig:EIem1} collects the results achieved by 4th-order compositions and constitutes the
analogous to Figure \ref{fig:em1}. We see that the overall efficiency of the schemes essentially corresponds to what is expected based on the effective
errors estimated in sections \ref{sect3} and \ref{sect4} and illustrated in the determination of traces for the linear problem \eqref{3.2a}. Figure \ref{fig:EIem2} is obtained by 6th-order methods. The results achieved by interchanging
the role of Euler explicit and Euler implicit are similar.

\begin{figure}
  \begin{center}
    \subfloat[]{
      \label{fig:EIem1a}
      \includegraphics[width=8cm]{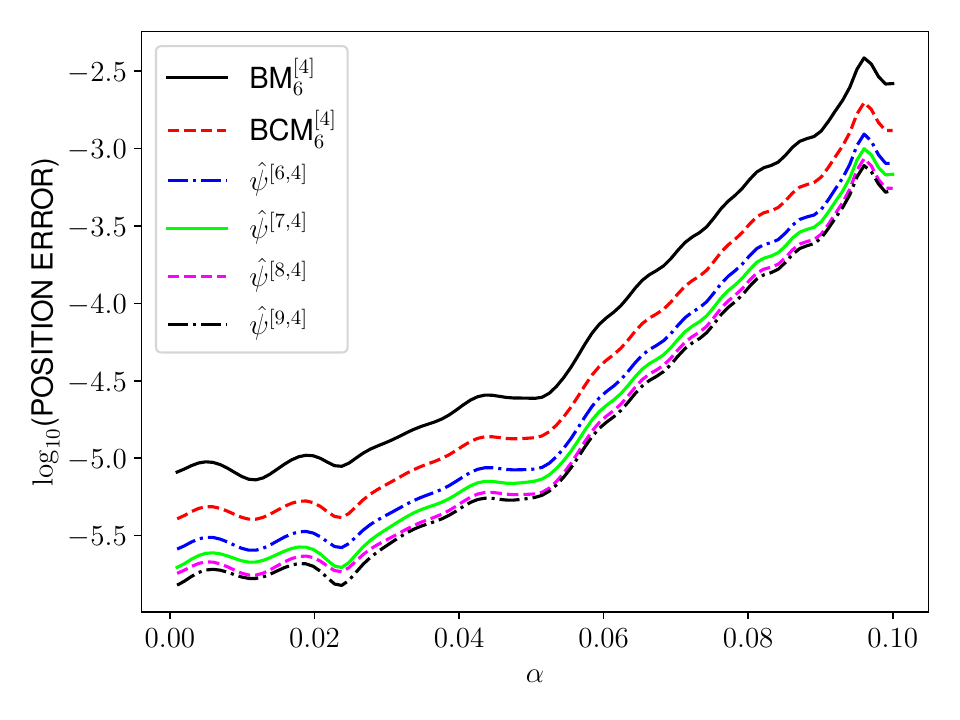}}
    \subfloat[]{
      \label{fig:EIem1b}
      \includegraphics[width=8cm]{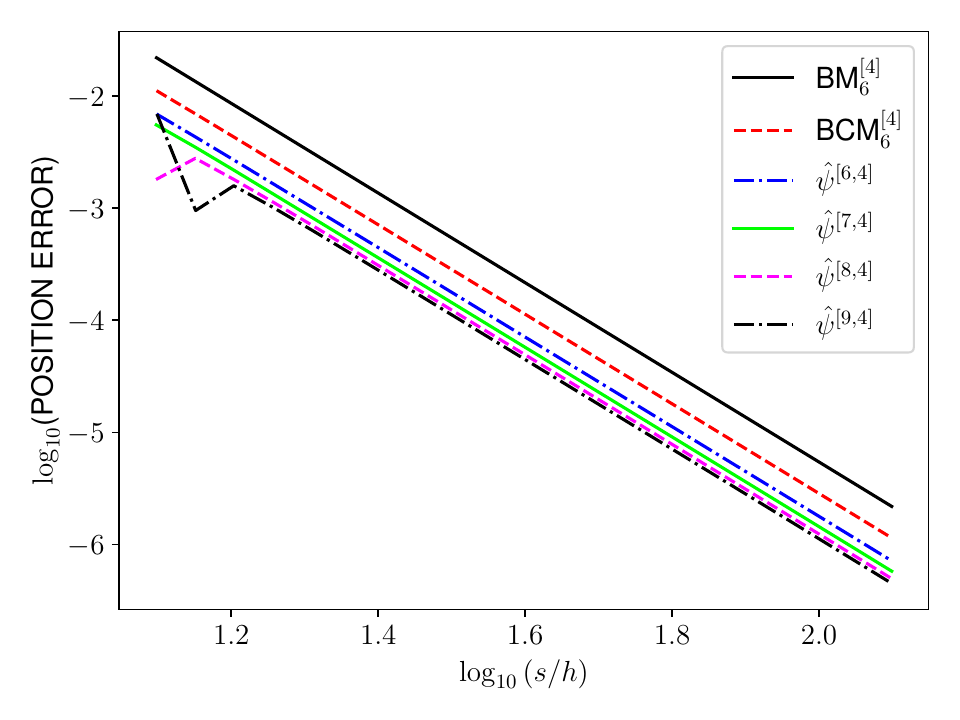}}
    \caption{\textbf{Motion of a charged particle under Lorentz force}. (a) Maximum position error for different values of the parameter $\alpha$ and for various 4th-order methods at the final time $t_f=200$, evaluated with the same computational cost $s/h=40$.
    $\chi_h$ and $\chi_h^*$ correspond to the explicit and implicit Euler methods, respectively.
      (b) Efficiency diagram for different 4th-order methods with $\alpha=0.07$ and $t_f=200$.}
    \label{fig:EIem1}
  \end{center}
\end{figure}

\begin{figure}
  \begin{center}
    \subfloat[]{
      \label{fig:EIem2a}
      \includegraphics[width=8cm]{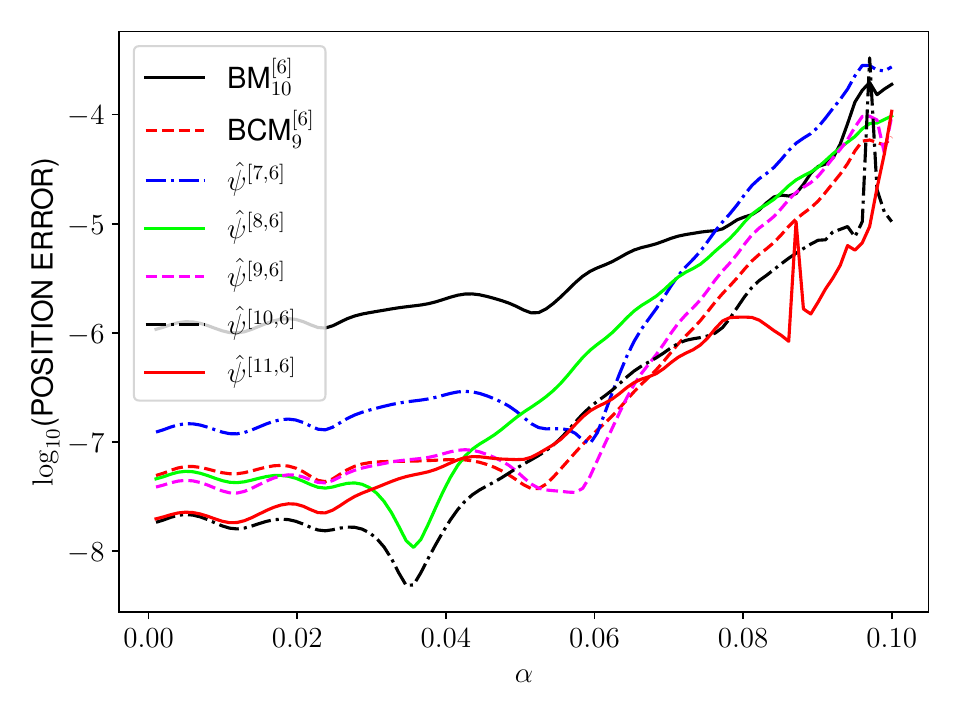}}
    \subfloat[]{
      \label{fig:EIem2b}
      \includegraphics[width=8cm]{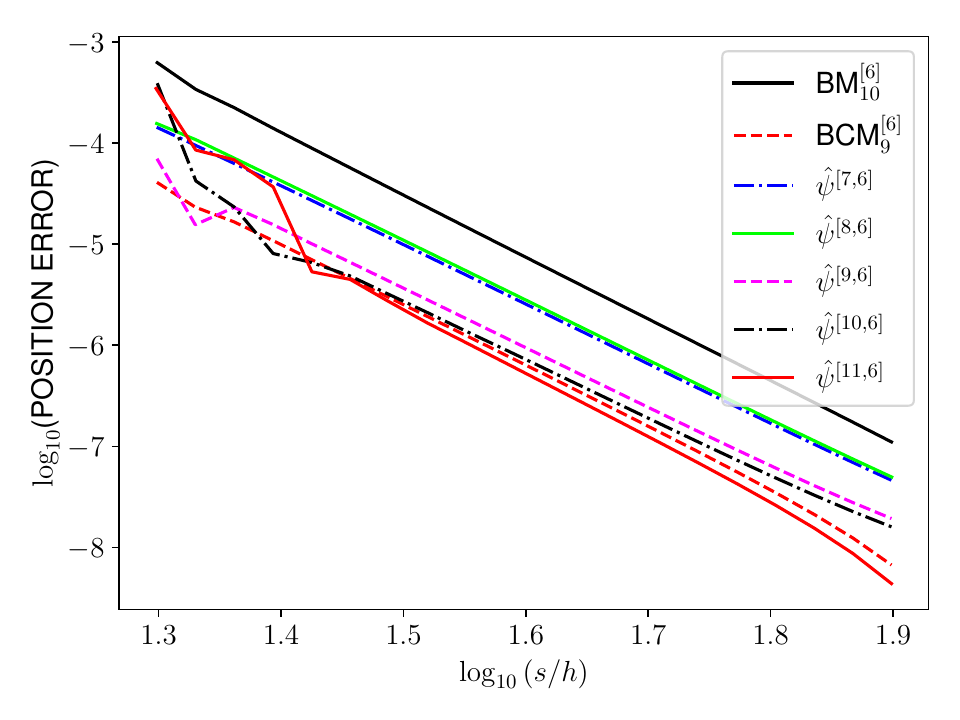}}
    \caption{\textbf{Motion of a charged particle under Lorentz force}. (a) Maximum position error for different values of the parameter $\alpha$ and for various 4th-order methods  at the final time $t_f=200$, evaluated with the same computational cost $s/h=40$.
    Here $\chi_h$ and $\chi_h^*$ correspond to the explicit and implicit Euler methods, respectively.
      (b) Efficiency diagram for different 4th-order methods with $\alpha=0.07$ and $t_f=200$.}
    \label{fig:EIem2}
  \end{center}
\end{figure}




\section{Concluding remarks}
We have presented new families of processed splitting methods of order 4 and 6 especially designed to be applied to problems which can be separated into three or more parts,
each of them being explicitly integrable. The construction strategy is as follows. First we determine the kernel, taking more stages than strictly required for solving the
order conditions, so that the free parameters are chosen so as to minimize not only the first term in the asymptotic expansion of the truncation error, but also higher order
terms, while keeping the size of the coefficients of the method reasonably small by keeping track of the 1-norm. The methods thus obtained are applied by computing the
trace of the solution of a linear system defined by three different random matrices. This simple test allows us not only to check the effective order but also to discard schemes
with large error constants. Second, for the most successful kernels we determine a particular processor $\pi_h$ also as a composition of elementary maps, whereas its adjoint
$\pi_h^*$ is taken as an approximation for the inverse $\pi_h^{-1}$. In this way the overall integrator is still time-symmetric, leading to good preservation properties.

As the results gathered in Tables \ref{tab.effic4} and \ref{effic_6} show, increasing the number of stages allows one to get kernels with smaller effective errors and coefficients,
the observed pattern closely following the rule of thumb formulated in \cite{mclachlan02foh}. Nevertheless, this efficiency pattern is not always followed when the methods
are applied in practice, as shown by the examples collected here. This is specially true when the basic map $\chi_h$ is formed as a composition of the exact solution of
each subproblem. On the contrary, if $\chi_h$ is a composition of first-order approximations, then the observed results agree nicely with the theoretical pattern. 

A possible
explanation for this behavior is related to the way the kernels are optimized. Specifically, in the optimization process carried out here, we have assumed that, for methods of order $r$, each Lie operator 
$E_{r+1,j}, \ j=1,2,\ldots$ in the basis of Table \ref{table.2} contributes equally to the error, but this not what always happens in practice. To see this point, let us consider
a system that is separable into just two parts, say $A$ and $B$, and the Hall basis for the corresponding free Lie algebra generated by $A$ and $B$ (similar
considerations follow for any other Hall--Viennot basis \cite{viennot78adl}), which we denote as $\tilde{E}_{i,j}$. Then,
\[
  \tilde E_{1,1}=A, \quad, \tilde E_{1,2}=B, \quad \tilde E_{2,1}=[A,B], \quad \tilde E_{3,1}=[A,[A,B]], \quad \tilde{E}_{3,2}=[B,[A,B]],
\]
etc. Now, if we assume that the contribution of each $\tilde E_{r+1,j}, \ j=1,2,\ldots$ to the error is similar, this is clearly not the case for the terms $E_{r+1,j}, \ j=1,2,\ldots$ in 
previous basis of Table \ref{table.2}: in fact, for the 4th- and 6th-order methods, the elements $E_{r+1,1}=Y_{r+1}, \ r=4,6$ provide the smallest contribution to the error.
On the other hand, if $\chi_h$ is taken as the explicit Euler method and the linear problem $\dot{U} = A U$ is considered, then 
\[
  \chi_h=(I+hA)=\exp\left(hA-\frac12h^2A^2+\frac13h^3A^3+\ldots\right)
\]
so that $E_{n,1}=Y_n=(-1)^{n+1}\frac1{n}A^n$ and $E_{n,j}=0, \ j>1$, since $[Y_j,Y_k]=0$. In consequence, the only surviving term at order $r+1$ is $E_{r+1,1}$. Thus,
for problems which are close to linear, when $\chi_h$ is taken as the explicit (or implicit) Euler method, we expect an important contribution from this term to the overall error.
Notice that is the case for the examples examined here.

\appendix

\section{Appendix}
For reader's convenience, we collect in Table \ref{table.processor} the coefficients for two particular processors:
\[
  \pi^{(9,4)} = \psiba_{\beta_7h}\circ  \psib_{\beta_6 h} \circ \chi_{\beta_5 h}^* \circ \chi_{\beta_4 h} \circ \chi_{\beta_3 h}^*
\circ \chi_{\beta_2 h} \circ \chi_{\beta_1 h}^*,
\]
and
\[  
   \pi^{(11,6)} = \psiba_{\beta_{23}h}\circ  \psib_{\beta_{22} h} \circ \cdots \circ \circ \chi_{\beta_2 h} \circ \chi_{\beta_1 h}^*,
\]    
corresponding to the kernels with $s=9$ and $s=11$ of effective order 4 and 6 of Tables \ref{table.4} and \ref{table.5}, respectively.

\begin{table}[t]
  \centering
    \renewcommand\arraystretch{1.1}
    \begin{tabular}{lll}
      \multicolumn{3}{c}{$\pi^{(9,4)}$}\\
      \hline 
      $\beta_1= -0.28566586026506785$ &\qquad $\beta_2= 0.015761586550701766$ \\
      $\beta_3= -0.04362530065430363$ & \qquad $\beta_4 =  -0.03618407560045836$ \\
      $\beta_5= 0.05244978481197771$ &   \qquad $\beta_6=  0.28558661670075497$ \\
      $\beta_7 = 0.011677248456395364$ \\
      \hline 
                                    & & \\      
      \multicolumn{3}{c}{$\pi^{(11,6)}$}\\
      \hline 
      $\beta_1= 0.2861698495034459$ &\qquad $\beta_2=  0.4134261834337682$ \\
      $\beta_3= 0.10540576774873363$ & \qquad $\beta_4 = -0.04664449698814812$ \\
      $\beta_5= 0.05672335497036459$ & \qquad $\beta_6 = 0.4990659695885505$ \\
      $\beta_7= -0.3426195751795226$ & \qquad $\beta_8 =0.3464936779661353$ \\
      $\beta_9= -0.23813674914660654$ & \qquad $\beta_{10} = 0.24491881441628852$ \\
      $\beta_{11}= -0.49669544275221306$ & \qquad $\beta_{12} =  -0.3122980257722082$ \\
       $\beta_{13}= 0.03146400131096136$ & \qquad $\beta_{14} =   -0.030063016455253767$ \\ 
      $\beta_{15}= 0.31240611169589994$ & \qquad $\beta_{16} =   -0.10319811497811636$ \\ 
      $\beta_{17}= -0.42098894976942247$ & \qquad $\beta_{18} = -0.2839790222445134$ \\ 
      $\beta_{19}= -0.039440980719714046$ & \qquad $\beta_{20} =  -0.020860135690795974$ \\ 
      $\beta_{21}= 0.05463728247473808$ & \qquad $\beta_{22} =  -0.16673300456832169 $ \\ 
     $\beta_{23}= 0.1509465011559501$ & \qquad  \\ 
      \hline 
     \end{tabular}
  \caption{\small{Coefficients of processors for methods of order 4 and 6. \label{table.processor}}}
\end{table}

\subsection*{Acknowledgments}
This work has been
funded by Ministerio de Ciencia e Innovaci\'on (Spain) through project PID2022-136585NB-C21, 
MCIN/AEI/10.13039/501100011033/FEDER, UE, and also by Generalitat Valenciana (Spain) through project CIAICO/2021/180.

\bibliographystyle{siam}


\begin{thebibliography}{10}

\bibitem{blanes21sps}
{\sc S.~Blanes, M.~P. Calvo, F.~Casas, and J.~M. Sanz-Serna}, {\em
  Symmetrically processed splitting integrators for enhanced {H}amiltonian
  {M}onte {C}arlo sampling}, SIAM J. Sci. Comput., 43 (2021), pp.~A3357--A3371.

\bibitem{blanes16aci}
{\sc S.~Blanes and F.~Casas}, {\em A {C}oncise {I}ntroduction to {G}eometric
  {N}umerical {I}ntegration}, {CRC} Press, 2016.

\bibitem{blanes23esm}
{\sc S.~Blanes, F.~Casas, C.~Gonz\'alez, and M.~Thalhammer}, {\em Efficient
  spliting methods based on modified potentials: numerical integration of
  linear parabolic problems and imaginary time propagation of the
  {S}chr\"odinger equation}, Commun. Comput. Phys., 33 (2023), pp.~937--961.

\bibitem{blanes04otn}
{\sc S.~Blanes, F.~Casas, and A.~Murua}, {\em On the numerical integration of
  ordinary differential equations by processed methods}, SIAM J. Numer. Anal.,
  42 (2004), pp.~531--552.

\bibitem{blanes06cmf}
{\sc S.~Blanes, F.~Casas, and A.~Murua}, {\em Composition methods
  for differential equations with processing}, SIAM J. Sci. Comput., 27 (2006),
  pp.~1817--1843.

\bibitem{blanes08sac}
{\sc S.~Blanes, F.~Casas, and A.~Murua}, {\em Splitting and
  composition methods in the numerical integration of differential equations},
  Bol. Soc. Esp. Mat. Apl., 45 (2008), pp.~89--145.

\bibitem{blanes24smf}
{\sc S.~Blanes, F.~Casas, and A.~Murua}, {\em Splitting methods
  for differential equations}, Acta Numerica, (in press) (2024).

\bibitem{blanes99siw}
{\sc S.~Blanes, F.~Casas, and J.~Ros}, {\em Symplectic integrators with
  processing: a general study}, SIAM J. Sci. Comput., 21 (1999), pp.~711--727.

\bibitem{blanes02psp}
{\sc S.~Blanes and P. C.~Moan}, {\em Practical symplectic partitioned
  {R}unge--{K}utta and {R}unge--{K}utta--{N}ystr\"om methods}, J. Comput. Appl.
  Math., 142 (2002), pp.~313--330.

\bibitem{butcher96tno}
{\sc J. C.~Butcher and J. M.~Sanz-Serna}, {\em The number of conditions for a
  {R}unge--{K}utta method to have effective order $p$}, Appl. Numer. Math., 22
  (1996), pp.~103--111.

\bibitem{hairer06gni}
{\sc E.~Hairer, C.~Lubich, and G.~Wanner}, {\em Geometric {N}umerical
  {I}ntegration. {S}tructure-{P}reserving {A}lgorithms for {O}rdinary
  {D}ifferential {E}quations}, Springer-Verlag, {S}econd~ed., 2006.

\bibitem{hairer93sod}
{\sc E.~Hairer, S.~N{\o}rsett, and G.~Wanner}, {\em Solving {O}rdinary
  {D}ifferential {E}quations {I}, {N}onstiff {P}roblems}, Springer-Verlag,
  {S}econd revised~ed., 1993.

\bibitem{he15vpa}
{\sc Y.~He, Y.~Sun, J.~Liu, and H.~Qin}, {\em Volume-preserving algorithms for
  charged particle dynamics}, J. Comput. Phys., 281 (2015), pp.~135--147.

\bibitem{he16hov}
{\sc Y.~He, Y.~Sun, J.~Liu, and H.~Qin}, {\em Higher order
  volume-preserving schemes for charged particle dynamics}, J. Comput. Phys.,
  305 (2016), pp.~172--184.

\bibitem{lopezmarcos97esi}
{\sc M. A.~L\'opez-Marcos, J. M.~Sanz-Serna, and R. D.~Skeel}, {\em Explicit symplectic
  integrators using {H}essian-vector products}, SIAM J. Sci. Comput., 18
  (1997), pp.~223--238.

\bibitem{mclachlan95otn}
{\sc R. I.~McLachlan}, {\em On the numerical integration of {ODE}'s by symmetric
  composition methods}, SIAM J. Sci. Comput., 16 (1995), pp.~151--168.

\bibitem{mclachlan02foh}
{\sc R. I.~McLachlan}, {\em Families of
  high-order composition methods}, Numer. Algor., 31 (2002), pp.~233--246.

\bibitem{mclachlan02sm}
{\sc R. I.~McLachlan and R.~Quispel}, {\em Splitting methods}, Acta Numerica, 11
  (2002), pp.~341--434.

\bibitem{mclachlan22tsi}
{\sc R.~I. McLachlan}, {\em Tuning symplectic integrators is easy and
  worthwhile}, Comm. Comp. Phys., 31 (2022), pp.~987--996.

\bibitem{munthe-kaas99cia}
{\sc H.~Munthe-Kaas and B.~Owren}, {\em Computations in a free {L}ie algebra},
  Phil. Trans. Royal Soc. A, 357 (1999), pp.~957--981.

\bibitem{sanz-serna94nhp}
{\sc J. M.~Sanz-Serna and M. P.~Calvo}, {\em Numerical {H}amiltonian {P}roblems},
  Chapman {\&} Hall, 1994.

\bibitem{suzuki90fdo}
{\sc M.~Suzuki}, {\em Fractal decomposition of exponential operators with
  applications to many-body theories and {M}onte {C}arlo simulations}, Phys.
  Lett. A, 146 (1990), pp.~319--323.

\bibitem{varadarajan84lgl}
{\sc V. S.~Varadarajan}, {\em Lie {G}roups, {L}ie {A}lgebras, and {T}heir
  {R}epresentations}, Springer-Verlag, 1984.

\bibitem{viennot78adl}
{\sc G.~Viennot}, {\em Alg\`ebres de {L}ie {L}ibres et {M}ono\"ides {L}ibres}, Lecture Notes in Mathematics 691, Springer-Verlag, 1978.

\bibitem{wang21coe1}
{\sc Y.~Wang, W.~Sun, F.~Liu, and X.~Wu}, {\em Construction of explicit
  symplectic integrators in general relativity. {I}. {S}chwarzschild black
  holes}, Astrophys. J., 907 (2021), p.~66.

\bibitem{wang21coe2}
{\sc Y.~Wang, W.~Sun, F.~Liu, and X.~Wu}, {\em Construction of
  explicit symplectic integrators in general relativity. {I}{I}.
  {R}eissner--{N}ordstr{\"o}m black holes}, Astrophys. J., 909 (2021), p.~22.

\bibitem{wisdom96sco}
{\sc J.~Wisdom, M.~Holman, and J.~Touma}, {\em Symplectic correctors}, in
  Integration Algorithms and Classical Mechanics, J.~Marsden, G.~Patrick, and
  W.~Shadwick, eds., vol.~10 of Fields Institute Communications, {A}merican
  {M}athematical {S}ociety, 1996, pp.~217--244.

\bibitem{yoshida90coh}
{\sc H.~Yoshida}, {\em Construction of higher order symplectic integrators},
  Phys. Lett. A, 150 (1990), pp.~262--268.

\end{thebibliography}

\end{document}